\documentclass[11pt,leqno]{amsart}
\usepackage[margin=20truemm, includefoot]{geometry}
\usepackage[utf8]{inputenc}
\usepackage{amsmath,amssymb,mathrsfs,stmaryrd}
\usepackage{bm}
\usepackage{enumerate}
\usepackage{amsthm}
\usepackage{hyperref}
\usepackage{comment}
\usepackage[all]{xy}
\usepackage{ascmac}
\usepackage{fancybox}
\usepackage{graphicx}
\usepackage{tikz-cd}
\usepackage{amsfonts}
\usepackage[bbgreekl]{mathbbol}
\usepackage{physics}
\usepackage{datetime}
\usepackage{shuffle}
\allowdisplaybreaks
\DeclareSymbolFontAlphabet{\mathbb}{AMSb}
\DeclareSymbolFontAlphabet{\mathbbl}{bbold}
\DeclareMathSymbol{\bbepsilon}{\mathord}{bbold}{"0F}
\xyoption{poly}
\newtheoremstyle{mystyle}
    {}
    {}
    {\normalfont}
    {}
    {\bfseries}
    {}
    { }
    {}

\DeclareMathOperator{\polylog}{Li}

\newtheorem{mainthm}{Theorem}
\newtheorem{defi}{Definition}[section]
\newtheorem{prop}[defi]{Proposition}

\newtheorem{cor}[defi]{Corollary}
\newtheorem{theo}[defi]{Theorem}

\newtheorem{lem}[defi]{Lemma}

\newtheorem{notation}[defi]{Notation}

\theoremstyle{remark}
\newtheorem{exa}[defi]{Example}
\newtheorem{rem}[defi]{Remarks}
\newtheorem*{rem*}{Remarks}

\makeatletter

\@addtoreset{equation}{section}
\makeatother

\def\C{{\mathbb C}}
\def\Z{{\mathbb Z}}

\def\Q{{\mathbb Q}}

\def\fH{\mathfrak{H}}

\title[Exponential generating functions of multiple polylogarithms]{Exponential generating functions of multiple polylogarithms at arbitrary integer indices}
\author{Kohei Kitamura}
\address[]{Department of Mathematics, The University of Osaka, Toyonaka, Osaka 560‐0043, Japan.}
\email{u347758b@ecs.osaka-u.ac.jp; kitamurakmath@gmail.com}
\begin{document}
\begin{abstract}
    We give an explicit formula expressing any finite product of depth-one polylogarithms at arbitrary integer indices as a $\mathbb{Q}$-linear combination of shuffle-type multivariable multiple polylogarithms. The proof starts with an exponential generating function. For any fixed positive multi-index $\mathbf{s}=(s_1,\ldots,s_r)$, the substitution $z_i\mapsto z_i\exp(x_1+\cdots+x_i)$ in the corresponding shuffle-type multiple polylogarithm yields the exponential generating function for the functions at indices $\mathbf{s}-\mathbf{k}$, where $\mathbf{k}$ ranges over $\mathbb{Z}_{\geq 0}^r$. We combine this identity with the shuffle product for iterated integrals at positive indices and compare coefficients to obtain the product formula. The appendix gives an alternative proof using the Euler differential operators $z_i\partial/\partial z_i$.
\end{abstract}
\maketitle
\tableofcontents
\section{Introduction}
    In this paper, we investigate multivariable multiple polylogarithms of $\shuffle$-type from geometric and combinatorial viewpoints, allowing their indices to have non-positive components. We work with multiple polylogarithms of $\shuffle$-type. Unless otherwise stated, MPL refers to a multiple polylogarithm of $\shuffle$-type. For an index $\mathbf{s}=(s_1,\ldots,s_r)\in\mathbb{Z}^r$, we define the MPLs of $\shuffle$-type
    \begin{align}
      \polylog^{\shuffle}_{\mathbf{s}}(z_1,\ldots,z_r)
      &:=
      \sum_{m_1>\cdots>m_r>0}
      \frac{z_1^{m_1-m_2}\cdots z_{r-1}^{m_{r-1}-m_r}z_r^{m_r}}
           {m_1^{s_1}\cdots m_r^{s_r}}
      \label{eq:intro-mpl-definition}\\
      &=
      \sum_{l_1,\ldots,l_r\geq 1}
      \frac{z_1^{l_1}\cdots z_r^{l_r}}
           {(l_1+\cdots+l_r)^{s_1}(l_2+\cdots+l_r)^{s_2}\cdots l_r^{s_r}}.
      \notag
    \end{align}
    These series define holomorphic functions on the unit polydisc $D^r=\{(z_1,\ldots,z_r)\in\mathbb{C}^r\mid |z_i|<1\text{ for all }i\}$. Moreover, if all components of the index $\mathbf{s}$ are positive, then $\polylog^\shuffle_\mathbf{s}$ has the following iterated integral representation,
    \begin{align*}
        \polylog^\shuffle_\mathbf{s}(z_1,\ldots,z_r)&=\int^1_0\qty(\frac{dt}{t})^{\circ (s_1-1)}\circ\qty(\frac{dt}{z_1^{-1}-t})\circ\cdots\circ\qty(\frac{dt}{t})^{\circ (s_r-1)}\circ\qty(\frac{dt}{z_r^{-1}-t})\\
        &=\underbrace{\int^1_0\frac{dt}{t}\left(\int_0^t\frac{dt}{t}\cdots\right.}_{s_1-1\text{ times}}\left(\int^t_0\right.\frac{dt}{z_1^{-1}-t}\cdots\left(\vphantom{\int^t_0}\right.\underbrace{\int^t_0\frac{dt}{t}\cdots}_{s_r-1\text{ times}}\left(\left.\int^t_0\frac{dt}{z_r^{-1}-t}\right)\cdots\right).
    \end{align*}
    \par
    For positive indices, the shuffle product formula has a geometric interpretation in terms of Chen's iterated integrals. We would like to extend the shuffle product formula to arbitrary integer indices. In \cite{Ki25}, the author derived the following product formula for one-variable MPLs with indices consisting entirely of non-positive integers (see \cite[Corollary 4.7]{Ki25} and also Corollary \ref{multi-variable NPMPL product} (2)). In this paper, we show the product formula for multi-variable and general indices which is induced from the usual shuffle relation.\par   
    Our starting point is the exponential generating function for MPLs all of whose index components are non-positive, established by Kaneko and Tsumura in \cite[Lemma 5.9]{KT18}. In our convention for the order of the indices, their formula reads
    \begin{equation}
      \sum_{\mathbf{k}\in\mathbb{Z}_{\geq 0}^r}\polylog^{\shuffle}_{-\mathbf{k}}(z_1,\ldots,z_r)\frac{\bm{x}^{\mathbf{k}}}{\mathbf{k}!}=\prod_{i=1}^r\frac{z_i \exp(x_1+\cdots+x_i)}{1-z_i \exp(x_1+\cdots+x_i)},
      \label{eq:intro-kt-generating-function}
    \end{equation}
    where $\bm{x}^{\mathbf{k}}=x_1^{k_1}\cdots x_r^{k_r}$ and
    $\mathbf{k}!=k_1!\cdots k_r!$.
    
    We first interpret this generating function as the composition of a multivariable MPL with exponential changes of its arguments. Indeed, the right-hand side of \eqref{eq:intro-kt-generating-function} can be written as
    \begin{equation*}
      \polylog^{\shuffle}_{\{0\}^r}
      (z_1\exp(x_1),z_2\exp(x_1+x_2),\ldots,z_r\exp(x_1+\cdots+x_r)).
    \end{equation*}
    We consider replacing $\{0\}^r$ with an arbitrary index $\mathbf{s}$ in this representation of multiple polylogarithms and expanding the resulting expression as a power series. We then obtain the following identity.
    \begin{mainthm}
        Let $\mathbf{s}\in\Z^r,\bm{x}\in\C^r$ and $z_1,\ldots,z_r\in\C$ such that each $|z_i|$ and $|z_i\exp(x_1+\cdots+x_i)|$ are sufficiently small. Then:
        \begin{equation}\label{eq:intro-shifted-generating-function}
            \polylog^{\shuffle}_{\mathbf{s}}(z_1\exp(x_1),z_2\exp(x_1+x_2),\ldots,z_r\exp(x_1+\cdots+x_r))=\sum_{\mathbf{k}\in\mathbb{Z}_{\geq 0}^r}\polylog^{\shuffle}_{\mathbf{s}-\mathbf{k}}(z_1,\ldots,z_r)\frac{\bm{x}^{\mathbf{k}}}{\mathbf{k}!}.
        \end{equation}
    \end{mainthm}
    \par
    Next, for $\mathbf{s}\in\mathbb{Z}_{>0}^r$, we apply the shuffle product formula to the product $\prod_{i=1}^r\polylog_{s_i}(z_ie^{x_i})$, that is
    \begin{align*}
        \prod_{i=1}^r\polylog_{s_i}(z_ie^{x_i})&=\prod_{i=1}^r\qty(\int^1_0\omega_0^{\circ (s_i-1)}\circ\omega_i)\\
        &=\int^1_0(\omega_0^{\circ (s_1-1)}\circ\omega_1)\shuffle\cdots\shuffle(\omega_0^{\circ (s_r-1)}\circ\omega_r),
    \end{align*}
    where $\omega_0:=dt/t$ and $\omega_i:=dt/(z_i^{-1}e^{-x_i}-t)$. Expanding both sides by the preceding generating series and comparing the coefficients of $\bm{x}^{\mathbf{j}}/\mathbf{j}!$, we obtain Theorem \ref{main theorem mixed shuffle}. To state the result, let $T_0,T_1,\ldots,T_r$ be noncommutative letters and define the shuffle coefficients $C_\sigma(\mathbf{s};\mathbf{n})$ by
    \begin{equation}
        (T_0^{s_1-1}T_1)\shuffle\cdots\shuffle(T_0^{s_r-1}T_r)=        \sum_{\sigma\in S_r}\sum_{\substack{\mathbf{n}\in\mathbb{Z}_{>0}^r                  |\mathbf{n}|=|\mathbf{s}|}}C_\sigma(\mathbf{s};\mathbf{n})T_0^{n_1-1}T_{\sigma(1)}\cdots T_0^{n_r-1}T_{\sigma(r)}.
      \label{eq:intro-shuffle-coefficients}
    \end{equation}
    For tuples, we use componentwise addition and subtraction, write $|\mathbf{a}|$ for the sum of their components, and denote the componentwise order by $\succeq$. We also set
    \begin{equation*}
      \binom{\mathbf{a}}{\mathbf{b}}
      :=\prod_{i=1}^r\binom{a_i}{b_i},
      \qquad
      \mathbf{j}_\sigma=(j_{\sigma(1)},\ldots,j_{\sigma(r)}),
      \qquad
      \bm{z}_\sigma=(z_{\sigma(1)},\ldots,z_{\sigma(r)}).
    \end{equation*}
    \begin{mainthm}\label{main theorem mixed shuffle}
    \textit{Let $\mathbf{s}=(s_1,\ldots,s_r)\in\mathbb{Z}_{>0}^r$ and $(z_1,\ldots,z_r)\in D^r$. Then, for every $\mathbf{j}=(j_1,\ldots,j_r)\in\mathbb{Z}_{\geq 0}^r$, one has}
    \begin{equation}\label{eq:intro-theorem-a}
        \prod_{i=1}^r\polylog_{s_i-j_i}(z_i)=\sum_{\sigma\in S_r}
        \sum_{\substack{\mathbf{n}\in\mathbb{Z}_{>0}^r\\|\mathbf{n}|=|\mathbf{s}|}}C_\sigma(\mathbf{s};\mathbf{n})\sum_{\substack{\mathbf{i}\in\mathbb{Z}_{\geq 0}^{r-1}\\\mathbf{j}_\sigma\succeq(\mathbf{i},0)}}(-1)^{|\mathbf{i}|}\binom{\mathbf{j}_\sigma}{(\mathbf{i},0)}\polylog^{\shuffle}_{(\mathbf{i},0)-(0,\mathbf{i})-\mathbf{j}_\sigma+\mathbf{n}}(\bm{z}_\sigma).
    \end{equation}
    Here $(\mathbf{i},0)$ and $(0,\mathbf{i})$ denote the tuples obtained by adjoining a zero at the end and at the beginning, respectively.
    \end{mainthm}
    
    Every integer tuple can be written as $\mathbf{s}-\mathbf{j}$ with $\mathbf{s}\in\mathbb{Z}_{>0}^r$ and $\mathbf{j}\in\mathbb{Z}_{\geq 0}^r$, so Theorem \ref{main theorem mixed shuffle} gives product formulae for depth-one polylogarithms at arbitrary integer indices. The resulting family of relations is geometric in the sense that it is derived from the shuffle product for iterated integrals. The generating series transfers these identities from positive indices to relations involving non-positive components, while retaining their origin in the geometry of iterated integration.
    
    Appendix A gives an alternative proof of Theorem \ref{main theorem mixed shuffle} using the commuting Euler differential operators
    \begin{equation*}
      D_i=z_i\frac{\partial}{\partial z_i},\qquad (i=1,\ldots,r).
    \end{equation*}
    Applying these operators to the shuffle identities at positive indices produces the same index shifts and binomial coefficients, and provides a differential approach to the relations obtained by coefficient comparison. From this approach, we can derive the formulae including the author's result:
    \begin{mainthm}\label{multi-variable NPMPL product}
        Let $\mathbf{j}=(j_1,\ldots,j_r)\in\Z_{\geq 0}^r$ and $z_1,\ldots,z_r\in\C$ satisfying $|z_i|<1$ for each $i$. Then we have:
        \begin{enumerate}[(1)]
            \item 
            \begin{equation}\label{eq.multi-variable NPMPL product}
                \prod_{i=1}^r\polylog_{-j_i}(z_i)=\sum_{\substack{\mathbf{j}\succeq(\mathbf{i},0),\mathbf{i}\in\Z^{r-1}_{\geq 0}}}(-1)^{\abs{\mathbf{i}}}\binom{\mathbf{j}}{(\mathbf{i},0)}\polylog^\shuffle_{(\mathbf{i},0)-(0,\mathbf{i})-\mathbf{j}}(z_1,\ldots,z_r),
            \end{equation}
            \item(\cite[Corollary 4.7]{Ki25})
            \begin{equation}\label{eq.mono-variable NPMPL product}
                \prod_{i=1}^r\polylog_{-j_i}(z)=\sum_{\substack{\mathbf{j}\succeq(\mathbf{i},0),\mathbf{i}\in\Z^{r-1}_{\geq 0}}}(-1)^{\abs{\mathbf{i}}}\binom{\mathbf{j}}{(\mathbf{i},0)}\polylog_{(\mathbf{i},0)-(0,\mathbf{i})-\mathbf{j}}(z).
            \end{equation}
        \end{enumerate}
    \end{mainthm}
    \begingroup
    \makeatletter
    \let\@tocwrite\@gobbletwo
    \makeatother
    
    \subsection*{Acknowledgements}
    The author would like to thank Koji Tasaka for fruitful discussions that inspired the ideas underlying this work, and Densuke Shiraishi for valuable comments that motivated this research. The author also would like to thank Hiroaki Nakamura for insightful suggestions and guidance throughout the preparation of this paper. This work was supported by JST SPRING, Grant Number JPMJSP2138.
    \subsection*{Declaration of Generative AI Use}
    During the preparation of this manuscript, the author used ChatGPT (OpenAI) to assist with English translation, language editing, and LaTeX code generation and editing. The author has reviewed and edited the generated content and takes full responsibility for the final manuscript.
    \par
    
    \endgroup
    \section{Generating series of MPLs at general indices}
    In this section, we will show the following proposition and focus on some corollaries.
    \begin{defi}[multi-variable MPLs, cf. (54) in \cite{KT18}]\label{MPL of each type}
        Let $s_1,\ldots,s_r\in\C$ and $z_1,\ldots,z_r\in\C$ such that $|z_i|<1$ for each $i=1,\ldots,r$. We define the \textup{multi-variable multiple polylogarithm} (for short MPL) 
        $\polylog^\shuffle_{s_1,\ldots,s_r}(z_1,\ldots,z_r)$ by
        \begin{align}
            \polylog_{s_1,\ldots,s_r}^\shuffle(z_1,\ldots,z_r):=&\sum_{m_1>\cdots>m_r>0}\frac{z_1^{m_1-m_2}\cdots z_{r-1}^{m_{r-1}-m_r}z_r^{m_r}}{m_1^{s_1}\cdots m_{r-1}^{s_{r-1}}m_r^{s_r}}\\
            =&\sum_{l_1,\ldots,l_r\geq1}\frac{z_1^{l_1}\cdots z_r^{l_r}}{(l_1+\cdots+l_r)^{s_1}(l_2+\cdots+l_r)^{s_2}\cdots l_r^{s_r}}.\notag
        \end{align}
    \end{defi}
    \begin{rem}
        (i) Kaneko--Tsumura (\cite{KT18}) defined another type of MPL $\polylog^\ast$ as 
        \begin{equation}
            \polylog^\ast_{s_1,\ldots,s_r}(z_1,\ldots,z_r):=\sum_{m_1>\cdots>m_r>0}\frac{z_1^{m_1}\cdots z_r^{m_r}}{m_1^{s_1}\cdots m_r^{s_r}}\qquad(|z_i|<1).
        \end{equation}
        They called the two types of MPLs $\polylog^\ast,\ \polylog^\shuffle$, of $\ast$-type, of $\shuffle$-type respectively.
        \par
        (ii) Clearly, one has
        \begin{equation}\label{one variable polylog}
            \polylog^\ast_{s_1,\ldots,s_r}(z,1,\ldots,1)=\polylog^\shuffle_{s_1,\ldots,s_r}(z,z,\cdots,z).
        \end{equation}
        We denote this common function by $\polylog_{s_1,\ldots,s_r}(z)$.
    \end{rem}
    \begin{prop}[Lemma 5.9 in \cite{KT18}]\label{generating function of NPMPL}
        Let $z_1,\ldots,z_r\in\C$ and $x_1,\ldots,x_r\in\C\ (r\in\Z_{>0})$. If $|z_i|\ (1\leq i\leq r)$ is sufficiently small, then the following holds:
        \begin{equation}\label{eq.generating function of NPMPL}
            \prod_{j=1}^r\frac{z_j\exp\qty(\sum_{\nu=1}^jx_\nu)}{1-z_j\exp\qty(\sum_{\nu=1}^jx_\nu)}=\sum_{k_1,\ldots,k_r\geq0}\polylog_{-k_1,\ldots,-k_r}^\shuffle(z_1,\ldots,z_r)\frac{x_1^{k_1}\cdots x_r^{k_r}}{k_1!\cdots k_r!}.
        \end{equation}
    \end{prop}
    \begin{proof}
        Refer to the proof of \cite[Lemma 5.9]{KT18} (note the order of indices and variables of MPL).
    \end{proof}
    In particular, setting $z_1=\cdots=z_r=z$ in \eqref{eq.generating function of NPMPL}, we obtain the generating function:
    \begin{equation}\label{eq.generating function of mono-variable NPMPL}
        \prod_{j=1}^r\frac{z\exp\qty(\sum_{\nu=1}^jx_\nu)}{1-z\exp\qty(\sum_{\nu=1}^jx_\nu)}=\sum_{k_1,\ldots,k_r\geq0}\polylog_{-k_1,\ldots,-k_r}^\shuffle(z)\frac{x_1^{k_1}\cdots x_r^{k_r}}{k_1!\cdots k_r!}.
    \end{equation}
    \begin{cor}
        Let $z\in\C$ and $x\in\C$. If $|z|$ is sufficiently small, then
        \begin{equation}\label{generating mono-polylog}
            \frac{z\exp(x)}{1-z\exp(x)}=\sum_{k\geq 0}\polylog_{-k}^\shuffle(z)\frac{x^k}{k!}.
        \end{equation}
    \end{cor}
    \begin{proof}
        To complete the proof, consider $r=1$ in \eqref{eq.generating function of mono-variable NPMPL}.
    \end{proof}
    \begin{cor}
        Let $z\in\C$ and $x_1,\ldots,x_r\in\C$. If $|z|$ is sufficiently small, then
        \begin{equation}
            \frac{z\exp(x_1+\cdots+x_r)}{1-z\exp(x_1+\cdots+x_r)}=\sum_{k_1,\ldots,k_r\geq 0}\polylog_{-(k_1+\cdots+k_r)}^\shuffle(z)\frac{x_1^{k_1}\cdots x_r^{k_r}}{k_1!\cdots k_r!}.
        \end{equation}
    \end{cor}
    \begin{proof}
        Substituting $x=x_1+\cdots+x_r$ in \eqref{generating mono-polylog}, we get
        \begin{align*}
            \frac{z\exp(x_1+\cdots+x_r)}{1-z\exp(x_1+\cdots+x_r)}&=\sum_{k\geq 0}\polylog^\shuffle_{-k}(z)\frac{(x_1+\cdots+x_r)^k}{k!}\\
            &=\sum_{k\geq 0}\polylog^\shuffle_{-k}(z)\sum_{k_1+\cdots+k_r=k}\binom{k}{k_1,\ldots,k_r}\frac{x_1^{k_1}\cdots x_r^{k_r}}{k!}\\
            &=\sum_{k\geq 0}\polylog^\shuffle_{-k}(z)\sum_{k_1+\cdots+k_r=k}\frac{x_1^{k_1}\cdots x_r^{k_r}}{k_1!\cdots k_r!}\\
            &=\sum_{k_1,\ldots,k_r\geq 0}\polylog^\shuffle_{-(k_1+\cdots+k_r)}(z)\frac{x_1^{k_1}\cdots x_r^{k_r}}{k_1!\cdots k_r!}.
        \end{align*}
    \end{proof}
    In the rest of this section, we aim to generalize \eqref{eq.generating function of NPMPL} and discuss convergence of the generalized one.\par
    First of all, we generalize \eqref{eq.generating function of NPMPL}. We focus on the equation appearing in the proof of \cite[Lemma 5.9]{KT18}:
    \begin{equation}\label{start of generalization}
        \sum_{k_1,\ldots,k_r=0}^\infty\polylog_{-k_1,\ldots,-k_r}^\shuffle(z_1,\ldots,z_r)\frac{x_1^{k_1}\cdots x_r^{k_r}}{k_1!\cdots k_r!}=\sum_{l_1,\ldots,l_r=1}^\infty\prod_{j=1}^r\qty(z_j\exp\qty(\sum_{\nu=1}^jx_\nu))^{l_j}.
    \end{equation}
    Here note that the notation is adjusted to maintain compatibility with the one in this paper. Then, we can interpret the right hand side of \eqref{start of generalization} as:
    \begin{equation*}
        \polylog^\shuffle_{\{0\}^r}\qty(z_1\exp(x_1),z_2\exp(x_1+x_2),\ldots,z_r\exp(x_1+\cdots+x_r)).
    \end{equation*}
    This interpretation is essential for our discussion, and we indeed reach the following formula.
    \begin{theo}\label{generating function of shuffle type}
        Let $s_1,\ldots,s_r\in\Z,x_1,\ldots,x_r\in\C$ and $z_1,\ldots,z_r\in\C$ such that each $|z_i|$ and $|z_i\exp(x_1+\cdots+x_i)|$ are sufficiently small. Then:
        \begin{align}
            \polylog^\shuffle_{s_1,\ldots,s_r}(z_1\exp(x_1),z_2\exp(x_1+x_2),\ldots,z_r\exp(x_1+\cdots+x_r))\\
            =\sum_{k_1,\ldots,k_r=0}^\infty\polylog^\shuffle_{-k_1+s_1,\ldots,-k_r+s_r}(z_1,\ldots,z_r)\frac{x_1^{k_1}\cdots x_r^{k_r}}{k_1!\cdots k_r!}.\notag
        \end{align}
    \end{theo}
    \begin{proof}
        To complete the proof, weight each term in the proof of \cite[Lemma 5.9]{KT18}, that is,
        \begin{align*}
            &\sum_{k_1,\ldots,k_r=0}^\infty\polylog^\shuffle_{-k_1+s_1,\ldots,-k_r+s_r}(z_1,\ldots,z_r)\frac{x_1^{k_1}\cdots x_r^{k_r}}{k_1!\cdots k_r!}\\
            &=\sum_{k_1,\ldots,k_r=0}^\infty\sum_{l_1,\ldots,l_r=1}^\infty\prod_{j=1}^r\frac{\qty(\qty(\sum_{\mu=j}^rl_\mu)x_j)^{k_j}}{k_j!}\cdot\frac{z_j^{l_j}}{\qty(\sum_{\mu=j}^rl_\mu)^{s_j}}\\
            &=\sum_{l_1,\ldots,l_r=1}^\infty\prod_{j=1}^r\frac{z_j^{l_j}}{\qty(\sum_{\mu=j}^rl_\mu)^{s_j}}\prod_{\mu=j}^r\exp(l_\mu x_j)\\
            &=\sum_{l_1,\ldots,l_r=1}^\infty\prod_{j=1}^r\frac{\qty(z_j\exp\qty(\sum_{\nu=1}^jx_\nu))^{l_j}}{\qty(\sum_{\mu=j}^rl_\mu)^{s_j}}\\
            &=\polylog^\shuffle_{s_1,\ldots,s_r}(z_1\exp(x_1),z_2\exp(x_1+x_2),\ldots,z_r\exp(x_1+\cdots+x_r)).
        \end{align*}
    \end{proof}
    It is natural to consider that a similar identity can hold in the case of $\ast$-type. Indeed, we obtain the following equation:
    \begin{theo}\label{generating function of stuffle type}
        Let $s_1,\ldots,s_r\in\Z,x_1,\ldots,x_r\in\C$ and $z_1,\ldots,z_r\in\C$ such that each $|z_i|$ and $|z_i\exp(x_i)|$ are sufficiently small. Then:
        \begin{align}
            \polylog^\ast_{s_1,\ldots,s_r}(z_1\exp(x_1),\ldots,z_r\exp(x_r))
            =\sum_{k_1,\ldots,k_r=0}^\infty\polylog^\ast_{-k_1+s_1,\ldots,-k_r+s_r}(z_1,\ldots,z_r)\frac{x_1^{k_1}\cdots x_r^{k_r}}{k_1!\cdots k_r!}.
        \end{align}
    \end{theo}
    \begin{proof}
    We obtain the desired equation by a similar argument. Indeed,
        \begin{align*}
            &\sum_{k_1,\ldots,k_r=0}^\infty\polylog^\ast_{-k_1+s_1,\ldots,-k_r+s_r}(z_1,\ldots,z_r)\frac{x_1^{k_1}\cdots x_r^{k_r}}{k_1!\cdots k_r!}\\
            &=\sum_{k_1,\ldots,k_r=0}^\infty\sum_{m_1>\cdots>m_r>0}^\infty\prod_{j=1}^r\frac{\qty(m_jx_j)^{k_j}}{k_j!}\cdot\frac{z_j^{m_j}}{m_j^{s_j}}\\
            &=\sum_{m_1>\cdots>m_r>0}^\infty\prod_{j=1}^r\frac{z_j^{m_j}}{m_j^{s_j}}\prod_{j=1}^r\exp(m_jx_j)\\
            &=\sum_{m_1>\cdots>m_r>0}^\infty\prod_{j=1}^r\frac{(z_j\exp(x_j))^{m_j}}{m_j^{s_j}}\\
            &=\polylog^\ast_{s_1,\ldots,s_r}(z_1\exp(x_1),\ldots,z_r\exp(x_r)).
        \end{align*}
    \end{proof}
\section{Shuffle relation and iterated integrals}
    Let $\Lambda$ be a set and let $\{X_\lambda\}_{\lambda\in\Lambda}$ be the corresponding set of letters. We denote by $\fH_\Lambda=\mathbb Q\langle X_\lambda\mid \lambda\in\Lambda\rangle$ the non-commutative polynomial algebra equipped with the shuffle product $\shuffle$, recursively defined by
    \begin{align}
        w\shuffle 1=1\shuffle w=w,\quad X_\alpha w_1\shuffle X_\beta w_2=X_\alpha(w_1\shuffle X_\beta w_2)+X_\beta(X_\alpha w_1\shuffle w_2).
    \end{align}
    Let $\Lambda=\{z\in\mathbb C\mid |z|<1\}=D.$
    Then we put
    \begin{equation}
        \fH_D^0=\Q\oplus\bigoplus_{z\in D\setminus\{0\}}\fH_D X_z.
    \end{equation}
    Thus ($\fH_D^0,\shuffle$) is a subalgebra of $(\fH_D,\shuffle)$.\par
    For $z\in D\setminus\{0\}$, set
    \begin{equation}
        \omega_0=\frac{dt}{t},\quad\omega_z=\frac{dt}{z^{-1}-t}.
    \end{equation}
    Since $z^{-1}\notin\overline{D}$, the form $\omega_z$ has no pole on the interval $[0,1]$. We use Chen's iterated integrals from $0$ to $1$. With this convention, for $s_1,\ldots,s_r\in\mathbb Z_{>0}$ and $z_1,\ldots,z_r\in D\setminus\{0\}$, one has
    \begin{align}
        \polylog^{\shuffle}_{s_1,\ldots,s_r}(z_1,\ldots,z_r)&=\int_{0}^{1}\omega_0^{\circ(s_1-1)}\circ\omega_{z_1}\circ\cdots\circ\omega_0^{\circ(s_r-1)}\circ\omega_{z_r}\\
        &=\underbrace{\int^1_0\frac{dt}{t}\left(\int_0^t\frac{dt}{t}\cdots\right.}_{s_1-1\text{ times}}\left(\int^t_0\right.\frac{dt}{z_1^{-1}-t}\cdots\left(\vphantom{\int^t_0}\right.\underbrace{\int^t_0\frac{dt}{t}\cdots}_{s_r-1\text{ times}}\left(\left.\int^t_0\frac{dt}{z_r^{-1}-t}\right)\cdots\right).\notag
    \end{align}
    Note that this iterated integral representation is well-defined for each $s_1,\ldots,s_r\in\Z_{>0}$. Equivalently, we define a $\mathbb Q$-linear map
    \begin{equation*}
        \polylog^{\shuffle}_{\bullet}:\fH_D^0\longrightarrow\C
    \end{equation*}
    by
    \begin{align}
        \begin{split}
            \polylog^{\shuffle}_{\bullet}(X_0^{s_1-1}X_{z_1}\cdots X_0^{s_r-1}X_{z_r}):=&\int_{0}^{1}\omega_0^{\circ(s_1-1)}\circ\omega_{z_1}\circ\cdots\circ\omega_0^{\circ(s_r-1)}\circ\omega_{z_r}\\
            =&\polylog^{\shuffle}_{s_1,\ldots,s_r}(z_1,\ldots,z_r),
        \end{split}
    \end{align}
    with the convention that $\polylog^\shuffle_\bullet(1)=1$.\par
    The following is the standard shuffle product formula for Chen's iterated
    integrals.
    \begin{theo}[Shuffle relation cf. \cite{Ch71}, \cite{Ho00}, \cite{Go01} etc.]\label{shuffle relation}
        The map
        \begin{equation*}
            \polylog^{\shuffle}_{\bullet}:\fH_D^0\longrightarrow\C
        \end{equation*}
        is a homomorphism of $\Q$-algebras. Namely, for any
        $u,v\in\fH_D^0$, one has
        \begin{equation}
            \polylog^{\shuffle}_{\bullet}(u\shuffle v)=\polylog^{\shuffle}_{\bullet}(u)\polylog^{\shuffle}_{\bullet}(v).
        \end{equation}
    \end{theo}
    \begin{proof}
        This follows from Chen's shuffle product formula for iterated integrals.
    \end{proof}
    \begin{rem}\label{analytic continuation of MPLs}
        The analytic continuation of MPLs is not needed in the proofs below. We use only the shuffle product formula in the domain $|z_i|<1$. The above iterated integral representation may be used to continue these functions, after choosing branches, to suitable simply connected subdomains of the complement of the singular divisors.
    \end{rem}
    Unless otherwise stated, the variables $z_i$ lie in the open unit disk $\mathbb D$. In proofs using iterated integrals, we first assume that all $z_i$ are nonzero; the identities extend to $\mathbb D^r$ by holomorphy.
\section{Products of two depth-one polylogarithms}
    In this section, we apply the shuffle relation to the case of depth $2$. Calculating
    \begin{align}\label{eq.words shuffle}
    \begin{split}
        X_0^{s_1-1}X_{z_1}\shuffle X_0^{s_2-1}X_{z_2}=&\sum_{i=0}^{s_2-1}\binom{s_1+i-1}{i}X_0^{s_1+i-1}X_{z_1}X_0^{s_2-i-1}X_{z_2}\\
        &+\sum_{j=0}^{s_1-1}\binom{s_2+j-1}{j}X_0^{s_2+j-1}X_{z_2}X_0^{s_1-j-1}X_{z_1},
    \end{split}
    \end{align}
    one has by shuffle relation (Theorem \ref{shuffle relation})
    \begin{equation}\label{shuffle of two depth 2}
    \begin{split}
        \polylog_{s_1}^\shuffle(z_1)\polylog_{s_2}^\shuffle(z_2)=&\sum_{i=0}^{s_2-1}\binom{s_1+i-1}{i}\polylog^\shuffle_{s_1+i,s_2-i}(z_1,z_2)\\
        &+\sum_{j=0}^{s_1-1}\binom{s_2+j-1}{j}\polylog^\shuffle_{s_2+j,s_1-j}(z_2,z_1).
    \end{split}
    \end{equation}
    \begin{lem}
        Let $z_1,z_2,x_1,x_2\in\C$ such that $|z_1|,|z_2|<1$, and $s_1,s_2\in\Z_{>0}$. Then for $(x_1,x_2)$ in a sufficiently small neighbourhood of $(0,0)$, the following identity holds
        \begin{equation}
            \polylog^\shuffle_{s_1,s_2}(z_1e^{x_1},z_2e^{x_2})=\sum_{j_1\geq k\geq 0,j_2\geq 0}(-1)^{k}\binom{j_1}{k}\cdot\polylog^\shuffle_{k-j_1+s_1,\,-k-j_2+s_2}(z_1,z_2)\frac{x_1^{j_1}x_2^{j_2}}{j_1!j_2!}.
        \end{equation}
    \end{lem}
    \begin{proof}
        From Theorem \ref{generating function of shuffle type} and \eqref{one variable polylog}, we verify
        \begin{align}
            \polylog^\shuffle_{s_1,s_2}(z_1e^{x_1},z_2e^{x_2})&=\polylog^\shuffle_{s_1,s_2}(z_1e^{x_1},z_2e^{x_1+(x_2-x_1)})\notag\\
            &=\sum_{i_1,i_2\geq 0}\polylog^\shuffle_{-i_1+s_1,-i_2+s_2}(z_1,z_2)\frac{x_1^{i_1}(x_2-x_1)^{i_2}}{i_1!i_2!}\notag\\
            &=\sum_{\substack{i_1,i_2\geq 0\\ i_2\geq j\geq 0}}\polylog^\shuffle_{-i_1+s_1,-i_2+s_2}(z_1,z_2)\cdot(-1)^{j}\binom{i_1+j}{j}\frac{x_1^{i_1+j}x_2^{i_2-j}}{(i_1+j)!(i_2-j)!}\notag\\
            &=\sum_{n\geq 0}\sum_{\substack{i_1+i_2=n\\i_1,i_2\geq 0}}\sum_{i_2\geq j\geq 0}\polylog^\shuffle_{-i_1+s_1,-i_2+s_2}(z_1,z_2)\cdot(-1)^{j}\binom{i_1+j}{j}\frac{x_1^{i_1+j}x_2^{i_2-j}}{(i_1+j)!(i_2-j)!}\notag\\
            &=\sum_{n\geq 0}\sum_{\substack{i_1+i_2=n\\i_1,i_2\geq 0}}\polylog^\shuffle_{-i_1+s_1,-i_2+s_2}(z_1,z_2)\sum_{n\geq j\geq i_1}(-1)^{j-i_1}\binom{j}{i_1}\frac{x_1^{j}x_2^{n-j}}{j!(n-j)!}\notag\\
            &=\sum_{\substack{n\geq i_1\geq 0}}\polylog^\shuffle_{-i_1+s_1,-n+i_1+s_2}(z_1,z_2)\sum_{n\geq j\geq i_1}(-1)^{j-i_1}\binom{j}{i_1}\frac{x_1^{j}x_2^{n-j}}{j!(n-j)!}\notag\\
            &=\sum_{n\geq j\geq i_1\geq 0}\polylog^\shuffle_{-i_1+s_1,-n+i_1+s_2}(z_1,z_2)\cdot(-1)^{j-i_1}\binom{j}{i_1}\frac{x_1^{j}x_2^{n-j}}{j!(n-j)!}\notag\\
            &=\sum_{j_1\geq k\geq 0,j_2\geq 0}(-1)^{j_1-k}\binom{j_1}{k}\cdot\polylog^\shuffle_{-k+s_1,\,k-j_1-j_2+s_2}(z_1,z_2)\frac{x_1^{j_1}x_2^{j_2}}{j_1!j_2!}\label{eq.before repracing}\\
            &=\sum_{j_1\geq k\geq 0,j_2\geq 0}(-1)^{k}\binom{j_1}{k}\cdot\polylog^\shuffle_{k-j_1+s_1,\,-k-j_2+s_2}(z_1,z_2)\frac{x_1^{j_1}x_2^{j_2}}{j_1!j_2!}.\label{eq.after repracing}
        \end{align}
    \end{proof}
    Here, replacing $k$ by $j_1-k$ in \eqref{eq.before repracing}, we obtain \eqref{eq.after repracing}, the main formula.
    \begin{theo}
        For any $j_1,j_2\in\Z_{\geq0}$ and $s_1,s_2\in\Z_{>0}$,
        \begin{equation}\label{eq.general index shuffle}
        \begin{split}
            \polylog_{-j_1+s_1}(z_1)\polylog_{-j_2+s_2}(z_2)=&\sum_{i=0}^{s_2-1}\sum_{k=0}^{j_1}(-1)^{k}\binom{s_1+i-1}{i}\binom{j_1}{k}\polylog^\shuffle_{k-j_1+s_1+i,\,-k-j_2+s_2-i}(z_1,z_2)\\
            &+\sum_{i=0}^{s_1-1}\sum_{k=0}^{j_2}(-1)^{k}\binom{s_2+i-1}{i}\binom{j_2}{k}\polylog^\shuffle_{k-j_2+s_2+i,\,-k-j_1+s_1-i}(z_2,z_1).
        \end{split}
        \end{equation}
    \end{theo}
    \begin{proof}
        Replacing $z_1,z_2$ by $z_1e^{x_1},z_2e^{x_2}$ respectively in \eqref{shuffle of two depth 2}, one has
        \begin{align}
        \begin{split}\label{calculating GF}
            \polylog_{s_1}(z_1e^{x_1})\polylog_{s_2}(z_2e^{x_2})=&\sum_{i=0}^{s_2-1}\binom{s_1+i-1}{i}\polylog^\shuffle_{s_1+i,s_2-i}(z_1e^{x_1},z_2e^{x_2})\\
            &+\sum_{j=0}^{s_1-1}\binom{s_2+j-1}{j}\polylog^\shuffle_{s_2+j,s_1-j}(z_2e^{x_2},z_1e^{x_1}).
        \end{split}
        \end{align}
        Expanding the left hand side in $x_1$ and $x_2$ by using Theorem \ref{generating function of shuffle type} gives
        \begin{equation*}
            \text{(LHS)}=\sum_{j_1,j_2\geq 0}\polylog_{-j_1+s_1}(z_1)\polylog_{-j_2+s_2}(z_2)\frac{x_1^{j_1}x_2^{j_2}}{j_1!j_2!}.
        \end{equation*}
        On the other hand, from the above lemma and Theorem \ref{generating function of shuffle type}, we obtain
        \begin{align}
        \begin{split}\label{s_1,s_2}
            &\sum_{i=0}^{s_2-1}\binom{s_1+i-1}{i}\polylog^\shuffle_{s_1+i,s_2-i}(z_1e^{x_1},z_2e^{x_2})\\
            &=\sum_{i=0}^{s_2-1}\binom{s_1+i-1}{i}\sum_{j_1\geq k\geq 0,j_2\geq 0}(-1)^{k}\binom{j_1}{k}\cdot\polylog^\shuffle_{k-j_1+s_1+i,\,-k-j_2+s_2-i}(z_1,z_2)\frac{x_1^{j_1}x_2^{j_2}}{j_1!j_2!},
        \end{split}\\
        \begin{split}\label{s_2,s_1}
            &\sum_{i=0}^{s_1-1}\binom{s_2+i-1}{i}\polylog^\shuffle_{s_2+i,s_1-i}(z_2e^{x_2},z_1e^{x_1})\\
            &=\sum_{i=0}^{s_1-1}\binom{s_2+i-1}{i}\sum_{j_2\geq k\geq 0,j_1\geq 0}(-1)^{k}\binom{j_2}{k}\cdot\polylog^\shuffle_{k-j_2+s_2+i,\,-k-j_1+s_1-i}(z_2,z_1)\frac{x_1^{j_1}x_2^{j_2}}{j_1!j_2!}.
        \end{split}
        \end{align}
        Since the right hand side of \eqref{calculating GF} is the sum of the above two terms \eqref{s_1,s_2} and \eqref{s_2,s_1}, comparing coefficients of $x_1^{j_1}x_2^{j_2}$ on both sides implies our assertion.
    \end{proof}
    \begin{exa}
        We consider the case $s_1-j_1=s_2-j_2=0$.\par
        (i) When $s_1=j_1=s_2=j_2=n>0$,
        \begin{align}
            \begin{split}
                \polylog_0(z_1)\polylog_0(z_2)=&\sum_{i=0}^{n-1}\sum_{k=0}^{n}(-1)^{k}\binom{n+i-1}{i}\binom{n}{k}\polylog^\shuffle_{k+i,\,-k-i}(z_1,z_2)\\
                &+\sum_{i=0}^{n-1}\sum_{k=0}^{n}(-1)^{k}\binom{n+i-1}{i}\binom{n}{k}\polylog^\shuffle_{k+i,\,-k-i}(z_2,z_1)
            \end{split}
            .
        \end{align}
        For instance, if $n=1$ then
        \begin{equation*}
            \polylog_0(z_1)\polylog_0(z_2)=\polylog^\shuffle_{0,0}(z_1,z_2)-\polylog^\shuffle_{1,-1}(z_1,z_2)+\polylog^\shuffle_{0,0}(z_2,z_1)-\polylog^\shuffle_{1,-1}(z_2,z_1),
        \end{equation*}
        and if $n=2$ then
        \begin{equation*}
        \begin{split}
            \polylog_0(z_1)\polylog_0(z_2)=&\polylog^\shuffle_{0,0}(z_1,z_2)-3\polylog^\shuffle_{2,-2}(z_1,z_2)+2\polylog^\shuffle_{3,-3}(z_1,z_2)\\
            &+\polylog^\shuffle_{0,0}(z_2,z_1)-3\polylog^\shuffle_{2,-2}(z_2,z_1)+2\polylog^\shuffle_{3,-3}(z_2,z_1)
        \end{split}
        \end{equation*}
        These equalities can be also verified from Faulhaber's formula.\par
        (ii) More generally, when $s_1=j_1=n>0,\,s_2=j_2=m>0,$
        \begin{align}
            \begin{split}
                \polylog_0(z_1)\polylog_0(z_2)=&\sum_{i=0}^{m-1}\sum_{k=0}^{n}(-1)^{k}\binom{n+i-1}{i}\binom{n}{k}\polylog^\shuffle_{k+i,\,-k-i}(z_1,z_2)\\
                &+\sum_{i=0}^{n-1}\sum_{k=0}^{m}(-1)^{k}\binom{m+i-1}{i}\binom{m}{k}\polylog^\shuffle_{k+i,\,-k-i}(z_2,z_1).
            \end{split}
        \end{align}
        For instance, if $n=1$ and $m=2$ then 
        \begin{align*}
        \begin{split}
            \polylog^\shuffle_{0}(z_1)\polylog^\shuffle_{0}(z_2)=&\polylog^\shuffle_{0,0}(z_1,z_2)-\polylog^\shuffle_{2,-2}(z_1,z_2)\\
            &+\polylog^\shuffle_{0,0}(z_2,z_1)-2\polylog^\shuffle_{1,-1}(z_2,z_1)+\polylog^\shuffle_{2,-2}(z_2,z_1).
        \end{split}
        \end{align*}
    \end{exa}
    \begin{exa}
        We consider the case $s_1-j_1=1,s_2-j_2=2$.\par
        (i) When $(s_1,j_1,s_2,j_2)=(1,0,2,0)$,
        \begin{equation}\label{ex.usuall shuffle}
        \begin{split}
            \polylog_{1}(z_1)\polylog_{2}(z_2)=&\sum_{i=0}^{1}\sum_{k=0}^{0}(-1)^{k}\binom{i}{i}\binom{0}{k}\polylog^\shuffle_{k+1+i,\,-k+2-i}(z_1,z_2)\\
            &+\sum_{i=0}^{0}\sum_{k=0}^{0}(-1)^{k}\binom{1+i}{i}\binom{0}{k}\polylog^\shuffle_{k+2+i,\,-k+1-i}(z_2,z_1)\\
            =&\polylog^\shuffle_{1,2}(z_1,z_2)+\polylog^\shuffle_{2,1}(z_1,z_2)+\polylog^\shuffle_{2,1}(z_2,z_1).
        \end{split}
        \end{equation}
        This is the shuffle relation in the usual sense.\par
        (ii) When $(s_1,j_1,s_2,j_2)=(2,1,2,0)$, 
        \begin{align*}
        \begin{split}
            \polylog_{1}(z_1)\polylog_{2}(z_2)=&\sum_{i=0}^{1}\sum_{k=0}^{1}(-1)^{k}\binom{1+i}{i}\binom{1}{k}\polylog^\shuffle_{k+1+i,\,-k+2-i}(z_1,z_2)\\
            &+\sum_{i=0}^{1}\sum_{k=0}^{0}(-1)^{k}\binom{1+i}{i}\binom{0}{k}\polylog^\shuffle_{k+2+i,\,-k+1-i}(z_2,z_1)
        \end{split}\\
        \begin{split}
            =&\polylog^\shuffle_{1,2}(z_1,z_2)+\polylog^\shuffle_{2,1}(z_1,z_2)-2\polylog^\shuffle_{3,0}(z_1,z_2)\\
            &+\polylog^\shuffle_{2,1}(z_2,z_1)+2\polylog^\shuffle_{3,0}(z_2,z_1).
        \end{split}
        \end{align*}\par
        (iii) When $(s_1,j_1,s_2,j_2)=(2,1,3,1)$,
        \begin{align}
        \begin{split}
            \polylog_{1}(z_1)\polylog_{2}(z_2)=&\sum_{i=0}^{2}\sum_{k=0}^{1}(-1)^{k}\binom{1+i}{i}\binom{1}{k}\polylog^\shuffle_{k+1+i,\,-k+2-i}(z_1,z_2)\\
            &+\sum_{i=0}^{1}\sum_{k=0}^{1}(-1)^{k}\binom{2+i}{i}\binom{1}{k}\polylog^\shuffle_{k+2+i,\,-k+1-i}(z_2,z_1)\notag
        \end{split}\\
        \begin{split}
            =&\polylog^\shuffle_{1,2}(z_1,z_2)+\polylog^\shuffle_{2,1}(z_1,z_2)+\polylog^\shuffle_{3,0}(z_1,z_2)-3\polylog^\shuffle_{4,-1}(z_1,z_2)\\
            &+\polylog^\shuffle_{2,1}(z_2,z_1)+2\polylog^\shuffle_{3,0}(z_2,z_1)-3\polylog^\shuffle_{4,-1}(z_2,z_1)
        \end{split}
        \end{align}
        The usual shuffle relation \eqref{ex.usuall shuffle} and this derive the following equation
        \begin{equation}\label{ex.2131}
            \polylog^\shuffle_{3,0}(z_1,z_2)=\polylog^\shuffle_{4,-1}(z_1,z_2)+\polylog^\shuffle_{4,-1}(z_2,z_1).
        \end{equation}\par
        (iv) When $(s_1,j_1,s_2,j_2)=(3,2,3,1)$,
        \begin{align*}
            \begin{split}
                \polylog_{1}(z_1)\polylog_{2}(z_2)=&\sum_{i=0}^{2}\sum_{k=0}^{2}(-1)^{k}\binom{2+i}{i}\binom{2}{k}\polylog^\shuffle_{k+1+i,\,-k+2-i}(z_1,z_2)\\
                &+\sum_{i=0}^{2}\sum_{k=0}^{1}(-1)^{k}\binom{2+i}{i}\binom{1}{k}\polylog^\shuffle_{k+2+i,\,-k+1-i}(z_2,z_1)
            \end{split}\\
            \begin{split}
                =&\polylog^\shuffle_{1,2}(z_1,z_2)+\polylog^\shuffle_{2,1}(z_1,z_2)+\polylog^\shuffle_{3,0}(z_1,z_2)-9\polylog^\shuffle_{4,-1}(z_1,z_2)+6\polylog^\shuffle_{5,-2}(z_1,z_2)\\
                &+\polylog^\shuffle_{2,1}(z_2,z_1)+2\polylog^\shuffle_{3,0}(z_2,z_1)+3\polylog^\shuffle_{4,-1}(z_2,z_1)-6\polylog^\shuffle_{5,-2}(z_2,z_1).
            \end{split}
        \end{align*}
    \end{exa}
    \begin{lem}\label{binomial lemma}
        Let $s,j\in\Z_{\geq 0}$ such that $s\geq j$. For any $n\in\Z_{\geq 0}$, we have
        \begin{equation}\label{eq.coefficients}
            \sum_{\substack{k,i\geq 0\\k+i=n}}(-1)^k\binom{s+i}{i}\binom{j}{k}=\binom{s-j+n}{n}.
        \end{equation}
        Here, if $k>j$, then $\binom{j}{k}=0$.
    \end{lem}
    \begin{proof}
        For any $i\geq0$, the integer $\binom{s+i}{i}$ is given as the coefficient of $x^i$ of the formal power series
        \begin{equation*}
            \frac{1}{(1-x)^{s+1}}=\sum_{i\geq0}\binom{s+i}{i}x^i.
        \end{equation*}
        Moreover, for any $0\leq k\leq j$, the integer $(-1)^k\binom{j}{k}$ appears as the coefficient of $x^k$ of the polynomial
        \begin{equation*}
            (1-x)^j=\sum_{k=0}^j(-1)^k\binom{j}{k}x^k.
        \end{equation*}
        Thus, we have
        \begin{equation*}
            \frac{1}{(1-x)^{s-j+1}}=\frac{(1-x)^j}{(1-x)^{s+1}}=\sum_{n\geq 0}\qty(\sum_{\substack{k,i\geq 0\\k+i=n}}(-1)^k\binom{s+i}{i}\binom{j}{k})x^n.
        \end{equation*}
        On the other hand, the coefficient of $x^n$ in $1/(1-x)^{s-j+1}$ is $\binom{s-j+n}{n}$ since $s-j\geq 0$ ;
        \begin{equation*}
            \frac{1}{(1-x)^{s-j+1}}=\sum_{n\geq 0}\binom{s-j+n}{n}x^n.
        \end{equation*}
        Hence we obtain the equation \eqref{eq.coefficients}.
    \end{proof}
    \begin{theo}
    
        Let $j_1,j_2\in\Z_{\geq 0}$ and $s_1,s_2\in\Z_{>0}$. Suppose that $s_1-j_1>0, s_2-j_2>0$. We define the sets 
        \begin{align*}
            R_1&:=\{(x,y)\in\Z^2\mid0\leq x\leq j_2,\,0\leq y\leq s_1-1\},\\
            R_2&:=\{(x,y)\in\Z^2\mid0\leq x\leq j_1,\,0\leq y\leq s_2-1\},\\
            A_i&:=\{(x,y)\in R_i\mid x+y\leq s_i-j_i-1\},\\
            B_i&:=R_i\setminus A_i.
        \end{align*}
        Then the following partial sum of the right hand side of \eqref{eq.general index shuffle}:
        \begin{equation}
        \begin{split}\label{eq.positive part}
            &\sum_{(k,i)\in A_2}(-1)^{k}\binom{s_1+i-1}{i}\binom{j_1}{k}\polylog^\shuffle_{k-j_1+s_1+i,\,-k-j_2+s_2-i}(z_1,z_2)\\
            &+\sum_{(k,i)\in A_1}(-1)^{k}\binom{s_2+i-1}{i}\binom{j_2}{k}\polylog^\shuffle_{k-j_2+s_2+i,\,-k-j_1+s_1-i}(z_2,z_1)
        \end{split}
        \end{equation}
        is equal to the linear sum derived from the shuffle product in the usual sense, that is,
        \begin{equation}
            \eqref{eq.positive part}=\polylog^\shuffle_\bullet(X_0^{s_1-j_1-1}X_{z_1}\shuffle X_0^{s_2-j_2-1}X_{z_2})(=\polylog^\shuffle_{s_1-j_1}(z_1)\polylog^\shuffle_{s_2-j_2}(z_2)).
        \end{equation}
        Consequently, we have
        \begin{equation}
            \begin{split}
                &\sum_{(k,i)\in B_2}(-1)^{k}\binom{s_1+i-1}{i}\binom{j_1}{k}\polylog^\shuffle_{k-j_1+s_1+i,\,-k-j_2+s_2-i}(z_1,z_2)\\
                &+\sum_{(k,i)\in B_1}(-1)^{k}\binom{s_2+i-1}{i}\binom{j_2}{k}\polylog^\shuffle_{k-j_2+s_2+i,\,-k-j_1+s_1-i}(z_2,z_1)
            \end{split}
            \qquad=0.
        \end{equation}
    \end{theo}
    \begin{proof}
        From Lemma \ref{binomial lemma}, \eqref{eq.words shuffle}, and \eqref{shuffle of two depth 2}, we have
        \begin{align*}
        \begin{split}
            &\sum_{\substack{0\leq k\leq j_1\\0\leq i\leq s_2-1\\k+i\leq s_2-j_2-1}}(-1)^{k}\binom{s_1+i-1}{i}\binom{j_1}{k}\polylog^\shuffle_{k-j_1+s_1+i,\,-k-j_2+s_2-i}(z_1,z_2)\\
            &+\sum_{\substack{0\leq k\leq j_2\\0\leq i\leq s_1-1\\k+i\leq s_1-j_1-1}}(-1)^{k}\binom{s_2+i-1}{i}\binom{j_2}{k}\polylog^\shuffle_{k-j_2+s_2+i,\,-k-j_1+s_1-i}(z_2,z_1)
        \end{split}\\
        \begin{split}
            =&\sum_{n=0}^{s_2-j_2-1}\sum_{\substack{k,i\geq 0\\k+i=n}}(-1)^{k}\binom{s_1+i-1}{i}\binom{j_1}{k}\polylog^\shuffle_{k-j_1+s_1+i,\,-k-j_2+s_2-i}(z_1,z_2)\\
            &+\sum_{n=0}^{s_1-j_1-1}\sum_{\substack{k,i\geq 0\\k+i=n}}(-1)^{k}\binom{s_2+i-1}{i}\binom{j_2}{k}\polylog^\shuffle_{k-j_2+s_2+i,\,-k-j_1+s_1-i}(z_2,z_1)
        \end{split}\\
        \begin{split}
            =&\sum_{n=0}^{s_2-j_2-1}\binom{s_1-j_1+n-1}{n}\polylog^\shuffle_{s_1-j_1+n,\,s_2-j_2-n}(z_1,z_2)\\
            &+\sum_{n=0}^{s_1-j_1-1}\binom{s_2-j_2+n-1}{n}\polylog^\shuffle_{s_2-j_2+n,\,s_1-j_1-n}(z_2,z_1)
        \end{split}\\
        =&\polylog^\shuffle_\bullet(X_0^{s_1-j_1-1}X_{z_1}\shuffle X_0^{s_2-j_2-1}X_{z_2}).
        \end{align*}
        This also implies that the latter assertion holds.
    \end{proof}
    \begin{exa}
        As a case that the differences $s_i-j_i$ are both non-positive, we consider $s_1-j_1=-1,\,s_2-j_2=-2$.
        \par
        (i) When $(s_1,j_1,s_2,j_2)=(1,2,1,3)$,
        \begin{align*}
            \begin{split}
                \polylog_{-1}(z_1)\polylog_{-2}(z_2)=&\sum_{i=0}^0\sum_{k=0}^2(-1)^k\binom{i}{i}\binom{2}{k}\polylog^\shuffle_{k-1+i,\,-k-2-i}(z_1,z_2)\\
                &+\sum_{i=0}^{0}\sum_{k=0}^{3}(-1)^{k}\binom{i}{i}\binom{3}{k}\polylog^\shuffle_{k-2+i,\,-k-1-i}(z_2,z_1)
            \end{split}
            \\
            \begin{split}
                =&\polylog^\shuffle_{-1,-2}(z_1,z_2)-2\polylog^\shuffle_{0,-3}(z_1,z_2)+\polylog^\shuffle_{1,-4}(z_1,z_2)\\
                &+\polylog^\shuffle_{-2,-1}(z_2,z_1)-3\polylog^\shuffle_{-1,-2}(z_2,z_1)+3\polylog^\shuffle_{0,-3}(z_2,z_1)-\polylog^\shuffle_{1,-4}(z_2,z_1).
            \end{split}
        \end{align*}
        \par
        (ii) When $(s_1,j_1,s_2,j_2)=(2,3,1,3)$,
        \begin{align*}
            \begin{split}
                \polylog_{-1}(z_1)\polylog_{-2}(z_2)=&\sum_{i=0}^0\sum_{k=0}^3(-1)^k\binom{i}{i}\binom{3}{k}\polylog^\shuffle_{k-1+i,\,-k-2-i}(z_1,z_2)\\
                &+\sum_{i=0}^{1}\sum_{k=0}^{3}(-1)^{k}\binom{1+i}{i}\binom{3}{k}\polylog^\shuffle_{k-2+i,\,-k-1-i}(z_2,z_1)
            \end{split}
            \\
            \begin{split}
                =&\polylog^\shuffle_{-1,-2}(z_1,z_2)-3\polylog^\shuffle_{0,-3}(z_1,z_2)+3\polylog^\shuffle_{1,-4}(z_1,z_2)-\polylog^\shuffle_{2,-5}(z_1,z_2)\\
                &+\polylog^\shuffle_{-2,-1}(z_2,z_1)-2\polylog^\shuffle_{-1,-2}(z_2,z_1)+2\polylog^\shuffle_{1,-4}(z_2,z_1)\\
                &\qquad-\polylog^\shuffle_{2,-5}(z_2,z_1)
            \end{split}
        \end{align*}\par
        (iii) When $(s_1,j_1,s_2,j_2)=(2,3,2,4)$, 
        \begin{align*}
            \begin{split}
                \polylog_{-1}(z_1)\polylog_{-2}(z_2)=&\sum_{i=0}^1\sum_{k=0}^3(-1)^k\binom{1+i}{i}\binom{3}{k}\polylog^\shuffle_{k-1+i,\,-k-2-i}(z_1,z_2)\\
                &+\sum_{i=0}^{1}\sum_{k=0}^{4}(-1)^{k}\binom{1+i}{i}\binom{4}{k}\polylog^\shuffle_{k-2+i,\,-k-1-i}(z_2,z_1)
            \end{split}
            \\
            \begin{split}
                =&\polylog^\shuffle_{-1,-2}(z_1,z_2)-\polylog^\shuffle_{0,-3}(z_1,z_2)-3\polylog^\shuffle_{1,-4}(z_1,z_2)+5\polylog^\shuffle_{2,-5}(z_1,z_2)\\
                &\qquad-2\polylog^\shuffle_{3,-6}(z_1,z_2)\\
                &+\polylog^\shuffle_{-2,-1}(z_2,z_1)-2\polylog^\shuffle_{-1,-2}(z_2,z_1)-2\polylog^\shuffle_{0,-3}(z_2,z_1)+8\polylog^\shuffle_{1,-4}(z_2,z_1)\\
                &\qquad-7\polylog^\shuffle_{2,-5}(z_2,z_1)+2\polylog^\shuffle_{3,-6}(z_2,z_1).
            \end{split}
        \end{align*}
        (iv) When $(s_1,j_1,s_2,j_2)=(3,4,2,4)$,
        \begin{align*}
        \begin{split}
            \polylog^\shuffle_{-1}(z_1)\polylog^\shuffle_{-2}(z_2)=&\sum_{i=0}^{1}\sum_{k=0}^{4}(-1)^k\binom{2+i}{i}\binom{4}{k}\polylog^\shuffle_{k-1+i,\,-k-2-i}(z_1,z_2)\\
            &+\sum_{i=0}^{2}\sum_{k=0}^{4}(-1)^{k}\binom{1+i}{i}\binom{4}{k}\polylog^\shuffle_{k-2+i,\,-k-1-i}(z_2,z_1)
        \end{split}
        \\
        \begin{split}
            =&\polylog^\shuffle_{-1,-2}(z_1,z_2)-\polylog^\shuffle_{0,-3}(z_1,z_2)-6\polylog^\shuffle_{1,-4}(z_1,z_2)+14\polylog^\shuffle_{2,-5}(z_1,z_2)\\
            &\qquad-11\polylog^\shuffle_{3,-6}(z_1,z_2)+3\polylog^\shuffle_{4,-7}(z_1,z_2)\\
            &+\polylog^\shuffle_{-2,-1}(z_2,z_1)-2\polylog^\shuffle_{-1,-2}(z_2,z_1)+\polylog^\shuffle_{0,-3}(z_2,z_1)-4\polylog^\shuffle_{1,-4}(z_2,z_1)\\
            &\qquad+11\polylog^\shuffle_{2,-5}(z_2,z_1)-10\polylog^\shuffle_{3,-6}(z_2,z_1)+3\polylog^\shuffle_{4,-7}(z_2,z_1).
        \end{split}
        \end{align*}
    \end{exa}
    \begin{exa}
        We consider the case $s_1-j_1=-1,s_2-j_2=2$, where exactly one of the differences $s_1-j_1$ and $s_2-j_2$ is positive, and the other is non-positive.\par
        (i) When $(s_1,j_1,s_2,j_2)=(1,2,3,1)$,
        \begin{align*}
        \begin{split}
            \polylog^\shuffle_{-1}(z_1)\polylog^\shuffle_{2}(z_2)=&\sum_{i=0}^{2}\sum_{k=0}^{2}(-1)^k\binom{i}{i}\binom{2}{k}\polylog^\shuffle_{k-1+i,\,-k+2-i}(z_1,z_2)\\
            &+\sum_{i=0}^{0}\sum_{k=0}^{1}(-1)^{k}\binom{2+i}{i}\binom{1}{k}\polylog^\shuffle_{k+2+i,\,-k-1-i}(z_2,z_1)
        \end{split}
        \\
        \begin{split}
            =&\polylog^\shuffle_{-1,2}(z_1,z_2)-\polylog^\shuffle_{0,1}(z_1,z_2)-\polylog^\shuffle_{2,-1}(z_1,z_2)+\polylog^\shuffle_{3,-2}(z_1,z_2)\\
            &+\polylog^\shuffle_{2,-1}(z_2,z_1)-\polylog^\shuffle_{3,-2}(z_2,z_1).
        \end{split}
        \end{align*}\par
        (ii) When $(s_1,j_1,s_2,j_2)=(2,3,3,1)$,
        \begin{align*}
        \begin{split}
            \polylog^\shuffle_{-1}(z_1)\polylog^\shuffle_{2}(z_2)=&\sum_{i=0}^{2}\sum_{k=0}^{3}(-1)^k\binom{1+i}{i}\binom{3}{k}\polylog^\shuffle_{k-1+i,\,-k+2-i}(z_1,z_2)\\
            &+\sum_{i=0}^{1}\sum_{k=0}^{1}(-1)^{k}\binom{2+i}{i}\binom{1}{k}\polylog^\shuffle_{k+2+i,\,-k-1-i}(z_2,z_1)
        \end{split}
        \\
        \begin{split}
            =&\polylog^\shuffle_{-1,2}(z_1,z_2)-\polylog^\shuffle_{0,1}(z_1,z_2)-4\polylog^\shuffle_{2,-1}(z_1,z_2)+7\polylog^\shuffle_{3,-2}(z_1,z_2)\\
            &\qquad-3\polylog^\shuffle_{4,-3}(z_1,z_2)\\
            &+\polylog^\shuffle_{2,-1}(z_2,z_1)+2\polylog^\shuffle_{3,-2}(z_2,z_1)-3\polylog^\shuffle_{4,-3}(z_2,z_1).
        \end{split}
        \end{align*}\par
        (iii) When $(s_1,j_1,s_2,j_2)=(2,3,4,2)$,
        \begin{align*}
        \begin{split}
            \polylog^\shuffle_{-1}(z_1)\polylog^\shuffle_{2}(z_2)=&\sum_{i=0}^{3}\sum_{k=0}^{3}(-1)^k\binom{1+i}{i}\binom{3}{k}\polylog^\shuffle_{k-1+i,\,-k+2-i}(z_1,z_2)\\
            &+\sum_{i=0}^{1}\sum_{k=0}^{2}(-1)^{k}\binom{3+i}{i}\binom{2}{k}\polylog^\shuffle_{k+2+i,\,-k-1-i}(z_2,z_1)
        \end{split}
        \\
        \begin{split}
            =&\polylog^\shuffle_{-1,2}(z_1,z_2)-\polylog^\shuffle_{0,1}(z_1,z_2)-5\polylog^\shuffle_{3,-2}(z_1,z_2)+9\polylog^\shuffle_{4,-3}(z_1,z_2)\\
            &\qquad-4\polylog^\shuffle_{5,-4}(z_1,z_2)\\
            &+\polylog^\shuffle_{2,-1}(z_2,z_1)+2\polylog^\shuffle_{3,-2}(z_2,z_1)-7\polylog^\shuffle_{4,-3}(z_2,z_1)+4\polylog^\shuffle_{5,-4}(z_2,z_1).
        \end{split}
        \end{align*}\par
        (iv) When $(s_1,j_1,s_2,j_2)=(3,4,4,2)$,
        \begin{align*}
        \begin{split}
            \polylog^\shuffle_{-1}(z_1)\polylog^\shuffle_{2}(z_2)=&\sum_{i=0}^{3}\sum_{k=0}^{4}(-1)^k\binom{2+i}{i}\binom{4}{k}\polylog^\shuffle_{k-1+i,\,-k+2-i}(z_1,z_2)\\
            &+\sum_{i=0}^{2}\sum_{k=0}^{2}(-1)^{k}\binom{3+i}{i}\binom{2}{k}\polylog^\shuffle_{k+2+i,\,-k-1-i}(z_2,z_1)
        \end{split}
        \\
        \begin{split}
            =&\polylog^\shuffle_{-1,2}(z_1,z_2)-\polylog^\shuffle_{0,1}(z_1,z_2)-15\polylog^\shuffle_{3,-2}(z_1,z_2)+39\polylog^\shuffle_{4,-3}(z_1,z_2)\\
            &\qquad-34\polylog^\shuffle_{5,-4}(z_1,z_2)+10\polylog^\shuffle_{6,-5}(z_1,z_2)\\
            &+\polylog^\shuffle_{2,-1}(z_2,z_1)+2\polylog^\shuffle_{3,-2}(z_2,z_1)+3\polylog^\shuffle_{4,-3}(z_2,z_1)-16\polylog^\shuffle_{5,-4}(z_2,z_1)\\
            &\qquad+10\polylog^\shuffle_{6,-5}(z_2,z_1).
        \end{split}
        \end{align*}
    \end{exa}
\section{Products of an arbitrary number of depth-one polylogarithms}
    \begin{notation}
        Let $r\in\Z_{>0}$, and let $\mathbf{i}=(i_1,\ldots,i_r)\in\Z_{\geq 0}^r,\ \mathbf{j}=(j_1,\ldots,j_r)\in\Z_{\geq 0}^r$. Then we write
        \begin{equation*}
            \mathbf{i}\pm\mathbf{j}:=(i_1\pm j_1,\ldots,i_r\pm j_r),\ \binom{\mathbf{j}}{\mathbf{i}}:=\binom{j_1}{i_1}\cdots\binom{j_r}{i_r},\ \mathbf{j}!=j_1!\cdots j_r!.
        \end{equation*}
        For $n\in\Z,\mathbf{s}=(s_1,\ldots,s_r)\in\Z^r$, we write
        \begin{equation*}
            (\mathbf{s},n):=(s_1,\ldots,s_r,n),\ (n,\mathbf{s}):=(n,s_1,\ldots,s_r)\in\Z^{r+1}.
        \end{equation*}
    \end{notation}
    We define the partial order $\succeq$ in $\Z^r$ as
    \begin{equation*}
        \mathbf{j}=(j_1,\ldots,j_r)\succeq\mathbf{i}=(i_1,\ldots,i_r)\Longleftrightarrow j_1\geq i_1,\ldots,j_r\geq i_r.
    \end{equation*}
    \begin{lem}\label{multiple generating function}
        Let $r>1$, $z_1,\ldots,z_r,x_1,\ldots,x_r\in\C$ such that $|z_i|<1$, and $\mathbf{s}=(s_1,\ldots,s_r)\in\Z_{>0}^r$. For $\bm{x}=(x_1,\ldots,x_r)\in\C^r$ in a sufficiently small neighbourhood of $\bm{0}$, the following holds
        \begin{equation}
            \begin{split}
                \polylog^\shuffle_\mathbf{s}(z_1e^{x_1},\ldots,z_re^{x_r})
                =\sum_{\substack{\mathbf{j}\succeq(\mathbf{i},0)\\\mathbf{j}\in\Z_{\geq 0}^r,\mathbf{i}\in\Z^{r-1}_{\geq 0}}}\polylog^\shuffle_{(\mathbf{i},0)-(0,\mathbf{i})-\mathbf{j}+\mathbf{s}}(\bm{z})\cdot(-1)^{\abs{\mathbf{i}}}\binom{\mathbf{j}}{(\mathbf{i},0)}\frac{\bm{x}^\mathbf{j}}{\mathbf{j}!}
            \end{split}
            ,
        \end{equation}
        where $\bm{x}^\mathbf{j}=x_1^{j_1}\cdots x_r^{j_r}$ and $\bm{z}=(z_1,\ldots,z_r)$.
    \end{lem}
    \begin{proof}
    From Theorem \ref{generating function of shuffle type}, we can verify the desired equation as follows:
        \begin{align*}
            &\polylog^\shuffle_{s_1,\ldots,s_r}(z_1e^{x_1},\ldots,z_re^{x_r})\\
            &=\polylog^\shuffle_{s_1,\ldots,s_r}\qty(z_1\exp(x_1),z_2\exp(x_1+(x_2-x_1)),\ldots,z_r\exp(x_1+\sum_{i=1}^{r-1}(x_{i+1}-x_i)))\\
            &=\sum_{k_1,\ldots,k_r=0}^\infty\polylog^\shuffle_{-k_1+s_1,\ldots,-k_r+s_r}(\bm{z})\frac{x_1^{k_1}(x_2-x_1)^{k_2}\cdots (x_r-x_{r-1})^{k_r}}{k_1!\cdots k_r!}\\
            &=\sum_{k_1,\ldots,k_r=0}^\infty\polylog^\shuffle_{-k_1+s_1,\ldots,-k_r+s_r}(\bm{z})\sum_{i_2,\ldots,i_r=0}^{k_2,\ldots,k_r}\prod_{l=2}^r(-1)^{i_l}\binom{k_l}{i_l}\frac{x_1^{k_1+i_2}x_2^{k_2-i_2+i_3}\cdots x_{r-1}^{k_{r-1}-i_{r-1}+i_r}x_r^{k_r-i_r}}{k_1!k_2!\cdots k_{r-1}!k_r!}\\
            \begin{split}
                &=\sum_{k_1,\ldots,k_r=0}^\infty\polylog^\shuffle_{-k_1+s_1,\ldots,-k_r+s_r}(\bm{z})\sum_{i_2,\ldots,i_r=0}^{k_2,\ldots,k_r}\binom{k_1+i_2}{i_2}\qty(\prod_{l=2}^{r-1}(-1)^{i_l}\binom{k_l-i_l+i_{l+1}}{i_{l+1}})\cdot(-1)^{i_r}\\
                &\hspace{6cm}\cdot\frac{x_1^{k_1+i_2}x_2^{k_2-i_2+i_3}\cdots x_{r-1}^{k_{r-1}-i_{r-1}+i_r}x_r^{k_r-i_r}}{(k_1+i_2)!(k_2-i_2+i_3)!\cdots(k_{r-1}-i_{r-1}+i_r)!(k_r-i_r)!}
            \end{split}\\
            \begin{split}
                &=\sum_{j_1\geq i_2\geq 0,\ldots,j_{r-1}\geq i_r\geq 0, j_r\geq 0}\polylog^\shuffle_{i_2-j_1+s_1,i_3-i_2-j_2+s_2,\ldots,i_r-i_{r-1}-j_{r-1}+s_{r-1},-i_r-j_r+s_r}(\bm{z})\\
                &\hspace{5cm}\cdot\qty(\prod_{l=1}^{r-1}(-1)^{i_{l+1}}\binom{j_l}{i_{l+1}})\frac{x_1^{j_1}\cdots x_r^{j_r}}{j_1!\cdots j_r!}
            \end{split}\\
            \begin{split}
                &=\sum_{j_1\geq i_1\geq 0,\ldots,j_{r-1}\geq i_{r-1}\geq 0, j_r\geq 0}\polylog^\shuffle_{i_1-j_1+s_1,i_2-i_1-j_2+s_2,\ldots,i_{r-1}-i_{r-2}-j_{r-1}+s_{r-1},-i_{r-1}-j_r+s_r}(\bm{z})\\
                &\hspace{5cm}\cdot\qty(\prod_{l=1}^{r-1}(-1)^{i_{l}}\binom{j_l}{i_{l}})\binom{j_r}{0}\frac{x_1^{j_1}\cdots x_r^{j_r}}{j_1!\cdots j_r!}.
            \end{split}
        \end{align*}
        Here, $\bm{z}=(z_1,\ldots,z_r)$.
    \end{proof}
    Let $\mathcal{R}:=\Q\langle T_0,T_1,\ldots,T_n,\ldots\rangle$ be the non-commutative ring over $\Q$ with respect to the concatenation, and let $W(\mathcal{R})$ be the set of all monomials in $\mathcal{R}$. For any $r>0,\,\mathbf{s}:=(s_1,\ldots,s_r),\,\mathbf{n}:=(n_1,\ldots,n_r)\in\Z_{>0}^r$, and $\sigma\in\mathfrak{S}_r$, an element of the symmetric group of degree $r$, we write $C_\sigma(\mathbf{s};\mathbf{n})$ for the coefficient of $T_0^{n_1-1}T_{\sigma(1)}\cdots T_0^{n_r-1}T_{\sigma(r)}$ in $T_0^{s_1-1}T_1\shuffle\cdots\shuffle T_0^{s_r-1}T_r$, that is, $C_\sigma(\mathbf{s};\mathbf{n})$ are determined by
    \begin{equation}\label{eq.shuffle coefficients}
        T_0^{s_1-1}T_1\shuffle\cdots\shuffle T_0^{s_r-1}T_r=\sum_{\substack{\sigma\in\mathfrak{S}_r\\\mathbf{n}=(n_1,\ldots,n_r)\in\Z_{>0}^r}}C_\sigma(\mathbf{s};\mathbf{n})T_0^{n_1-1}T_{\sigma(1)}\cdots T_0^{n_r-1}T_{\sigma(r)}.
    \end{equation}
    Note that $C_\sigma(\mathbf{s};\mathbf{n})=0$ if $|\mathbf{n}|\neq|\mathbf{s}|$.
    \begin{theo}\label{mixed shuffle product}
        Let $\mathbf{s}=(s_1,\ldots,s_r)\in\Z^r_{>0}$, and $z_1,\ldots,z_r\in\C$ be complex variables such that $|z_i|<1$.
        Then for any $\mathbf{j}=(j_1,\ldots,j_r)\in\Z_{\geq0}^r$, we have
        \begin{equation}
            \prod_{i=1}^r\polylog_{s_i-j_i}(z_i)=\sum_{\sigma\in\mathfrak{S}_r}\sum_{\substack{\mathbf{n}\in\Z_{>0}^r\\|\mathbf{n}|=|\mathbf{s}|}}C_\sigma(\mathbf{s};\mathbf{n})\qty(\sum_{\substack{\mathbf{j}_\sigma\succeq(\mathbf{i},0)\\\mathbf{i}\in\Z^{r-1}_{\geq 0}}}(-1)^{\abs{\mathbf{i}}}\binom{\mathbf{j}_\sigma}{(\mathbf{i},0)}\polylog^\shuffle_{(\mathbf{i},0)-(0,\mathbf{i})-\mathbf{j}_\sigma+\mathbf{n}}(\bm{z}_\sigma)).
        \end{equation}
        where $\mathbf{j}_\sigma=(j_{\sigma(1)},\ldots,j_{\sigma(r)}),\,\bm{z}_\sigma=(z_{\sigma(1)},\ldots,z_{\sigma(r)})$. 
    \end{theo}
    \begin{proof}
        The case $r=1$ is immediate, so we assume $r\geq 2$. Firstly, in \eqref{eq.shuffle coefficients} we put $T_0=X_0,\,T_i=X_{z_ie^{x_i}}$ and verify the equation
        \begin{equation*}
            X_0^{s_1-1}X_{z_1e^{x_1}}\shuffle\cdots\shuffle X_0^{s_r-1}X_{z_re^{x_r}}=\sum_{\substack{\sigma\in\mathfrak{S}_r\\\mathbf{n}=(n_1,\ldots,n_r)\in\Z_{>0}^r}}C_\sigma(\mathbf{s};\mathbf{n})X_0^{n_1-1}X_{z_{\sigma(1)}e^{x_{\sigma(1)}}}\cdots X_0^{n_r-1}X_{z_{\sigma(r)}e^{x_{\sigma(r)}}}.
        \end{equation*}
        Thus, we can compute by using the usual shuffle relation (Theorem \ref{shuffle relation}) as:
        \begin{align*}
            \prod_{i=1}^r\polylog_{s_i}(z_ie^{x_i})&=\prod_{i=1}^r\polylog^\shuffle_\bullet(X_0^{s_i-1}X_{z_ie^{x_i}})\\
            &=\polylog^\shuffle_\bullet(X_0^{s_1-1}X_{z_1e^{x_1}}\shuffle\cdots\shuffle X_0^{s_r-1}X_{z_re^{x_r}})\\
            \begin{split}
                &=\sum_{\sigma\in\mathfrak{S}_r}\sum_{\substack{\mathbf{n}=(n_1,\ldots,n_r)\in\Z_{>0}^r\\|\mathbf{n}|=|\mathbf{s}|}}C_\sigma(\mathbf{s};\mathbf{n})\polylog^\shuffle_\bullet(X_0^{n_1-1}X_{z_{\sigma(1)}e^{x_{\sigma(1)}}}\cdots X_0^{n_r-1}X_{z_{\sigma(r)}e^{x_{\sigma(r)}}})
            \end{split}
            \\
            \begin{split}
                &=\sum_{\sigma\in\mathfrak{S}_r}\sum_{\substack{\mathbf{n}=(n_1,\ldots,n_r)\in\Z_{>0}^r\\|\mathbf{n}|=|\mathbf{s}|}}C_\sigma(\mathbf{s};\mathbf{n})\polylog^\shuffle_{n_1,\ldots,n_r}(z_{\sigma(1)}e^{x_{\sigma(1)}},\ldots,z_{\sigma(r)}e^{x_{\sigma(r)}}).
            \end{split}
        \end{align*}
        Here, note that repetitions are allowed in the expansion of the above formula. Namely, even if $X_{z_ie^{x_i}} = X_{z_je^{x_j}}$ for some $i$ and $j$, we formally regard $z_ie^{x_i}$ and $z_je^{x_j}$ as distinct and take the sum accordingly. Then, from Lemma \ref{multiple generating function}, we can verify
        \begin{align}
            \begin{split}
                &\sum_{\sigma\in\mathfrak{S}_r}\sum_{\substack{\mathbf{n}=(n_1,\ldots,n_r)\in\Z_{>0}^r\\|\mathbf{n}|=|\mathbf{s}|}}C_\sigma(\mathbf{s};\mathbf{n})\polylog^\shuffle_{n_1,\ldots,n_r}(z_{\sigma(1)}e^{x_{\sigma(1)}},\ldots,z_{\sigma(r)}e^{x_{\sigma(r)}})
            \end{split}\notag
            \\
            \begin{split}
                &=\sum_{\sigma\in\mathfrak{S}_r}\sum_{\substack{\mathbf{n}=(n_1,\ldots,n_r)\in\Z_{>0}^r\\|\mathbf{n}|=|\mathbf{s}|}}C_\sigma(\mathbf{s};\mathbf{n})\sum_{\substack{\mathbf{j}\succeq(\mathbf{i},0)\\\mathbf{j}\in\Z_{\geq 0}^r,\mathbf{i}\in\Z^{r-1}_{\geq 0}}}\polylog^\shuffle_{(\mathbf{i},0)-(0,\mathbf{i})-\mathbf{j}+\mathbf{n}}(\bm{z}_\sigma)
                \cdot(-1)^{\abs{\mathbf{i}}}\binom{\mathbf{j}}{(\mathbf{i},0)}\frac{\bm{x}_\sigma^{\mathbf{j}}}{\mathbf{j}!}.
            \end{split}\notag
            \\\label{eq.permutation(1)}
            \begin{split}
                &=\sum_{\sigma\in\mathfrak{S}_r}\sum_{\substack{\mathbf{n}\in\Z_{>0}^r\\|\mathbf{n}|=|\mathbf{s}|}}C_\sigma(\mathbf{s};\mathbf{n})\sum_{\substack{\mathbf{j}_\sigma\succeq(\mathbf{i},0)\\\mathbf{j}\in\Z_{\geq 0}^r,\mathbf{i}\in\Z^{r-1}_{\geq 0}}}\polylog^\shuffle_{(\mathbf{i},0)-(0,\mathbf{i})-\mathbf{j}_\sigma+\mathbf{n}}(\bm{z}_\sigma)
                \cdot(-1)^{\abs{\mathbf{i}}}\binom{\mathbf{j}_\sigma}{(\mathbf{i},0)}\frac{\bm{x}_\sigma^{\mathbf{j}_\sigma}}{\mathbf{j}_\sigma!}
            \end{split}
            \\\label{eq.permutation(2)}
            \begin{split}
                &=\sum_{\sigma\in\mathfrak{S}_r}\sum_{\substack{\mathbf{n}\in\Z_{>0}^r\\|\mathbf{n}|=|\mathbf{s}|}}C_\sigma(\mathbf{s};\mathbf{n})\sum_{\substack{\mathbf{j}_\sigma\succeq(\mathbf{i},0)\\\mathbf{j}\in\Z_{\geq 0}^r,\mathbf{i}\in\Z^{r-1}_{\geq 0}}}\polylog^\shuffle_{(\mathbf{i},0)-(0,\mathbf{i})-\mathbf{j}_\sigma+\mathbf{n}}(\bm{z}_\sigma)
                \cdot(-1)^{\abs{\mathbf{i}}}\binom{\mathbf{j}_\sigma}{(\mathbf{i},0)}\frac{\bm{x}^{\mathbf{j}}}{\mathbf{j}!}
            \end{split}
            \\
            \begin{split}
                &=\sum_{\sigma\in\mathfrak{S}_r}\sum_{\substack{\mathbf{n}\in\Z_{>0}^r\\|\mathbf{n}|=|\mathbf{s}|}}C_\sigma(\mathbf{s};\mathbf{n})\sum_{\mathbf{j}\in\Z_{\geq 0}^r}\qty(\sum_{\substack{\mathbf{j}_\sigma\succeq(\mathbf{i},0)\\\mathbf{i}\in\Z^{r-1}_{\geq 0}}}\polylog^\shuffle_{(\mathbf{i},0)-(0,\mathbf{i})-\mathbf{j}_\sigma+\mathbf{n}}(\bm{z}_\sigma)\cdot(-1)^{\abs{\mathbf{i}}}\binom{\mathbf{j}_\sigma}{(\mathbf{i},0)})\frac{\bm{x}^{\mathbf{j}}}{\mathbf{j}!}
            \end{split}\notag
            \\
            \begin{split}
                &=\sum_{\mathbf{j}\in\Z_{\geq 0}^r}\sum_{\sigma\in\mathfrak{S}_r}\sum_{\substack{\mathbf{n}\in\Z_{>0}^r\\|\mathbf{n}|=|\mathbf{s}|}}C_\sigma(\mathbf{s};\mathbf{n})\qty(\sum_{\substack{\mathbf{j}_\sigma\succeq(\mathbf{i},0)\\\mathbf{i}\in\Z^{r-1}_{\geq 0}}}\polylog^\shuffle_{(\mathbf{i},0)-(0,\mathbf{i})-\mathbf{j}_\sigma+\mathbf{n}}(\bm{z}_\sigma)\cdot(-1)^{\abs{\mathbf{i}}}\binom{\mathbf{j}_\sigma}{(\mathbf{i},0)})\frac{\bm{x}^{\mathbf{j}}}{\mathbf{j}!}
            \end{split}\notag
        \end{align}
        where $\bm{x}_\sigma^{\mathbf{j}}:=x_{\sigma(1)}^{j_{1}}\cdots x_{\sigma(r)}^{j_r}$, and the equality \eqref{eq.permutation(1)}$=$\eqref{eq.permutation(2)} is justified by
        \begin{equation*}
            \frac{\bm{x}_\sigma^{\mathbf{j}_\sigma}}{\mathbf{j}_\sigma!}=\frac{x_{\sigma(1)}^{j_{\sigma(1)}}\cdots x_{\sigma(r)}^{j_{\sigma(r)}}}{j_{\sigma(1)}!\cdots j_{\sigma(r)}!}=\frac{x_1^{j_1}\cdots x_r^{j_r}}{j_1!\cdots j_r!}=\frac{\bm{x}^\mathbf{j}}{\mathbf{j}!}.
        \end{equation*}
        \par
        On the other hand, one has clearly
        \begin{equation*}
            \prod_{i=1}^r\polylog^\shuffle_{s_i}(z_ie^{x_i})=\sum_{\mathbf{j}=(j_1,\ldots,j_r)\in\Z_{\geq 0}^r}\qty(\prod_{i=1}^r\polylog^\shuffle_{s_i-j_i}(z_i))\frac{\bm{x}^\mathbf{j}}{\mathbf{j}!}.
        \end{equation*}
        Hence we obtain the desired equation by comparing the coefficients of $\bm{x}^\mathbf{j}/\mathbf{j}!$ in the series.
    \end{proof}
    \begin{cor}
        Let $\mathbf{k}=(k_1,\ldots,k_r)\in\Z_{\geq 0}^r$ and $z_1,\ldots,z_r\in\C$ satisfying $|z_i|<1$ for each $i$. Then
        \begin{equation}\label{mixed shuffle at non-positive}
            \prod_{i=1}^r\polylog^\shuffle_{-k_i}(z_i)=\sum_{\sigma\in\mathfrak{S}_r}\sum_{\substack{\mathbf{k}_\sigma+(\{1\}^r)\succeq(\mathbf{i},0)\\\mathbf{i}\in\Z^{r-1}_{\geq 0}}}(-1)^{\abs{\mathbf{i}}}\binom{\mathbf{k}_\sigma+(\{1\}^r)}{(\mathbf{i},0)}\polylog^\shuffle_{(\mathbf{i},0)-(0,\mathbf{i})-\mathbf{k}_\sigma}(\bm{z}_\sigma),
        \end{equation}
        where $\mathbf{k}_\sigma=(k_{\sigma(1)},\ldots,k_{\sigma(r)})$. Moreover, if $z_1=\cdots=z_r=:z$, then
        \begin{equation}
            \prod_{i=1}^r\polylog^\shuffle_{-k_i}(z)=\sum_{\sigma\in\mathfrak{S}_r}\sum_{\substack{\mathbf{k}_\sigma+(\{1\}^r)\succeq(\mathbf{i},0)\\\mathbf{i}\in\Z^{r-1}_{\geq 0}}}(-1)^{\abs{\mathbf{i}}}\binom{\mathbf{k}_\sigma+(\{1\}^r)}{(\mathbf{i},0)}\polylog^\shuffle_{(\mathbf{i},0)-(0,\mathbf{i})-\mathbf{k}_\sigma}(\bm{z}),
        \end{equation}
        where we put $\bm{z}=(z,\ldots,z)$.
    \end{cor}
    \begin{proof}
        In Theorem \ref{mixed shuffle product}, put $\mathbf{s}=(\{1\}^r)$ and $\mathbf{j}=(k_1+1,\ldots,k_r+1)=\mathbf{k}+(\{1\}^r)$. Then $C_\sigma(\mathbf{s};\mathbf{n})=1$ holds for any $\sigma\in\mathfrak{S}_r$ and $\mathbf{n}=(\{1\}^r)$. Since $\{(\{1\}^r)\}=\{\mathbf{n}\in\Z_{>0}^r\mid|\mathbf{n}|=|\mathbf{s}|(=r)\}$, we have \eqref{mixed shuffle at non-positive}.
    \end{proof}
    \appendix
    \section{Another proof of Theorem \ref{main theorem mixed shuffle}}
    In this section, we shall show Theorem \ref{main theorem mixed shuffle} by using (partial) differential operators. This proof is more intuitive and simpler than the preceding one.\par
    For a holomorphic function $f(z_1,\ldots,z_r)$ on an open subset $U$ in $\C^r$, we define the partial differential operators
    \begin{equation}
        D_i:=z_i\frac{\partial}{\partial z_i}\quad(i=1,\ldots,r).
    \end{equation}
    These differential operators $D_i$ have $\C$-linear action on $\mathcal{O}(U)$ for each open subset $U\subset\C^r$, where $\mathcal{O}(U)$ is the ring of holomorphic functions on $U$.
    \begin{prop}\label{differential operator}
        \begin{enumerate}[(i)]
            \item 
            For $i,j=1,\ldots,r,\,D_i\circ D_j=D_j\circ D_i$,
            \item 
            For $s_1,\ldots,s_r\in\Z$ and $i=1,\ldots,r-1$,
            \begin{align}
                D_i(\polylog^\shuffle_{s_1,\ldots,s_r}(z_1,\ldots,z_r))&=\polylog^\shuffle_{s_1,\ldots,s_{i-1},s_i-1,s_{i+1},\ldots,s_r}(z_1,\ldots,z_r)-\polylog^\shuffle_{s_1,\ldots,s_i,s_{i+1}-1,s_{i+2},\ldots,s_r}(z_1,\ldots,z_r),\\
                D_r(\polylog^\shuffle_{s_1,\ldots,s_r}(z_1,\ldots,z_r))&=\polylog^\shuffle_{s_1,\ldots,s_r-1}(z_1,\ldots,z_r).
            \end{align}
        \end{enumerate}
    \end{prop}
    \begin{proof}
        $\polylog_{s_1,\ldots,s_r}^\shuffle(z_1,\ldots,z_r)$ is holomorphic on the unit polydisk $D^r$. Any holomorphic function is term by term differentiable on its series expansion. Hence we have
        \begin{align*}
            D_i(\polylog^\shuffle_{s_1,\ldots,s_r}(z_1,\ldots,z_r))&=D_i\left(\sum_{m_1>\cdots>m_r>0}\frac{z_1^{m_1-m_2}\cdots z_{r-1}^{m_{r-1}-m_r}z_r^{m_r}}{m_1^{s_1}\cdots m_{r-1}^{s_{r-1}}m_r^{s_r}}\right)\\
            &=\sum_{m_1>\cdots>m_r>0}(m_i-m_{i+1})\frac{z_1^{m_1-m_2}\cdots z_{r-1}^{m_{r-1}-m_r}z_r^{m_r}}{m_1^{s_1}\cdots m_{r-1}^{s_{r-1}}m_r^{s_r}}\\
            &=\sum_{m_1>\cdots>m_r>0}\frac{z_1^{m_1-m_2}\cdots z_{r-1}^{m_{r-1}-m_r}z_r^{m_r}}{m_1^{s_1}\cdots m_{i-1}^{s_{i-1}}m_i^{s_i-1}m_{i+1}^{s_{i+1}}\cdots m_r^{s_r}}\\
            &\qquad-\sum_{m_1>\cdots>m_r>0}\frac{z_1^{m_1-m_2}\cdots z_{r-1}^{m_{r-1}-m_r}z_r^{m_r}}{m_1^{s_1}\cdots m_{i}^{s_{i}}m_{i+1}^{s_{i+1}-1}m_{i+2}^{s_{i+2}}\cdots m_r^{s_r}}\\
            &=\polylog^\shuffle_{s_1,\ldots,s_{i-1},s_i-1,s_{i+1},\ldots,s_r}(z_1,\ldots,z_r)-\polylog^\shuffle_{s_1,\ldots,s_i,s_{i+1}-1,s_{i+2},\ldots,s_r}(z_1,\ldots,z_r).
        \end{align*}
        Similarly, we obtain
        \begin{align*}
            D_r(\polylog^\shuffle_{s_1,\ldots,s_r}(z_1,\ldots,z_r))&=D_r\left(\sum_{m_1>\cdots>m_r>0}\frac{z_1^{m_1-m_2}\cdots z_{r-1}^{m_{r-1}-m_r}z_r^{m_r}}{m_1^{s_1}\cdots m_{r-1}^{s_{r-1}}m_r^{s_r}}\right)\\
            &=\sum_{m_1>\cdots>m_r>0}m_r\frac{z_1^{m_1-m_2}\cdots z_{r-1}^{m_{r-1}-m_r}z_r^{m_r}}{m_1^{s_1}\cdots m_{r-1}^{s_{r-1}}m_r^{s_r}}\\
            &=\sum_{m_1>\cdots>m_r>0}\frac{z_1^{m_1-m_2}\cdots z_{r-1}^{m_{r-1}-m_r}z_r^{m_r}}{m_1^{s_1}\cdots m_{r-1}^{s_{r-1}}m_r^{s_r-1}}\\
            &=\polylog^\shuffle_{s_1,\ldots,s_r-1}(z_1,\ldots,z_r).
        \end{align*}
    \end{proof}
    \begin{cor}\label{index descending}
        For $\mathbf{j}=(j_1,\ldots,j_r)\in\Z_{\geq 0}^r$ and $\mathbf{s}=(s_1,\ldots,s_r)\in\Z^r$,
        \begin{equation}
            \left(\prod_{i=1}^rD_i^{\circ j_i}\right)(\polylog^\shuffle_\mathbf{s}(z_1,\ldots,z_r))=\sum_{\substack{\mathbf{j}\succeq(\mathbf{i},0),\mathbf{i}\in\Z^{r-1}_{\geq 0}}}(-1)^{\abs{\mathbf{i}}}\binom{\mathbf{j}}{(\mathbf{i},0)}\polylog^\shuffle_{(\mathbf{i},0)-(0,\mathbf{i})-\mathbf{j}+\mathbf{s}}(z_1,\ldots,z_r).
        \end{equation}
    \end{cor}
    \begin{proof}
        The case $r=1$ is immediate, so we assume $r\geq 2$. When $i=1,\ldots,r-1$, from Proposition \ref{differential operator} (ii), we have
        \begin{align*}
            D_i^{\circ j_i}(\polylog^\shuffle_{s_1,\ldots,s_r}(z_1,\ldots,z_r))=\sum_{k_i=0}^{j_i}(-1)^{k_i}\binom{j_i}{k_i}\polylog^\shuffle_{s_1,\ldots,s_i+k_i-j_i,\,s_{i+1}-k_i,\ldots,s_r}(z_1,\ldots,z_r).
        \end{align*}
        Thus, we recursively obtain
        \begin{align*}
            &\left(\prod_{i=1}^{r-1}D_i^{\circ j_i}\right)(\polylog^\shuffle_{s_1,\ldots,s_r}(z_1,\ldots,z_r))\\
            &=(D_1^{\circ j_1}\circ\cdots\circ D_{r-1}^{\circ j_{r-1}})(\polylog^\shuffle_{s_1,\ldots,s_r}(z_1,\ldots,z_r))\\
            &=(D_1^{\circ j_1}\circ\cdots\circ D_{r-2}^{\circ j_{r-2}})\qty(\sum_{k_{r-1}=0}^{j_{r-1}}(-1)^{k_{r-1}}\binom{j_{r-1}}{k_{r-1}}\polylog^\shuffle_{s_1,\ldots,s_{r-1}+k_{r-1}-j_{r-1},\,s_r-k_{r-1}}(z_1,\ldots,z_r))\\
            &\quad\cdots\notag\\
            \begin{split}
                &=\sum_{j_1\geq k_1\geq 0,\ldots,j_{r-1}\geq k_{r-1}\geq 0}(-1)^{k_1+\cdots+k_{r-1}}\binom{j_1}{k_1}\cdots\binom{j_{r-1}}{k_{r-1}}\\
                &\hspace{5cm}\cdot\polylog^\shuffle_{s_1+k_1-j_1,s_2+k_2-k_1-j_2,\cdots,s_{r-1}+k_{r-1}-k_{r-2}-j_{r-1},s_r-k_{r-1}}(z_1,\ldots,z_r).
            \end{split}
        \end{align*}
        Applying Proposition \ref{differential operator} once again, we have
        \begin{align*}
            &\left(\prod_{i=1}^{r}D_i^{\circ j_i}\right)(\polylog^\shuffle_{s_1,\ldots,s_r}(z_1,\ldots,z_r))\\
            \begin{split}
                &=D_r^{\circ j_r}\left(\sum_{j_1\geq k_1\geq 0,\ldots,j_{r-1}\geq k_{r-1}\geq 0}(-1)^{k_1+\cdots+k_{r-1}}\binom{j_1}{k_1}\cdots\binom{j_{r-1}}{k_{r-1}}\right.\\
                &\hspace{5cm}\left.\phantom{\sum_J}\cdot\polylog^\shuffle_{s_1+k_1-j_1,s_2+k_2-k_1-j_2,\cdots,s_{r-1}+k_{r-1}-k_{r-2}-j_{r-1},s_r-k_{r-1}}(z_1,\ldots,z_r)\right)
            \end{split}\\
            \begin{split}
                &=\sum_{j_1\geq k_1\geq 0,\ldots,j_{r-1}\geq k_{r-1}\geq 0}(-1)^{k_1+\cdots+k_{r-1}}\binom{j_1}{k_1}\cdots\binom{j_{r-1}}{k_{r-1}}\binom{j_r}{0}\\
                &\hspace{5cm}\cdot\polylog^\shuffle_{s_1+k_1-j_1,s_2+k_2-k_1-j_2,\cdots,s_{r-1}+k_{r-1}-k_{r-2}-j_{r-1},s_r-k_{r-1}-j_r}(z_1,\ldots,z_r)
            \end{split}\\
            &=\sum_{\substack{\mathbf{j}\succeq(\mathbf{i},0),\mathbf{i}\in\Z^{r-1}_{\geq 0}}}
                (-1)^{\abs{\mathbf{i}}}\binom{\mathbf{j}}{(\mathbf{i},0)}\polylog^\shuffle_{(\mathbf{i},0)-(0,\mathbf{i})-\mathbf{j}+\mathbf{s}}(z_1,\ldots,z_r).
        \end{align*}
    \end{proof}
    \begin{rem}
        Note that the components of the indices $\mathbf{s}$ in this corollary may take arbitrary integer values, not only positive integer values.
    \end{rem}
    \begin{proof}[Proof of Theorem \ref{main theorem mixed shuffle} (Theorem \ref{mixed shuffle product})]
        First of all, for $i,k=1,\ldots,r$ we have
        \begin{equation*}
            D_i^{\circ j_i}(\polylog_{s_k}(z_k))=\left\lbrace
            \begin{array}{ll}
                \delta_{ik}\polylog_{s_k-j_i}(z_k)&(\text{if }j_i>0),\\
                \polylog_{s_k}(z_k)&(\text{if } j_i=0),
            \end{array}
            \right.
        \end{equation*}
        where $\delta_{ik}$ is Kronecker's delta. Thus we obtain
        \begin{align}\label{eq.product of derivatives}
            \left(\prod_{i=1}^rD_i^{\circ j_i}\right)\left(\prod_{i=1}^r\polylog_{s_i}(z_i)\right)=\prod_{i=1}^r\polylog_{s_i-j_i}(z_i).
        \end{align}
        On the other hand, the usual shuffle product (Theorem \ref{shuffle relation}), together with \eqref{eq.shuffle coefficients}, gives
        \begin{align}
            \prod_{i=1}^r\polylog_{s_i}(z_i)&=\prod_{i=1}^r\polylog^\shuffle_\bullet(X_0^{s_i-1}X_{z_i})\label{eq.product(1)}\\
            &=\polylog^\shuffle_\bullet(X_0^{s_1-1}X_{z_1}\shuffle\cdots\shuffle X_0^{s_r-1}X_{z_r})\notag\\
            &=\sum_{\sigma\in\mathfrak{S}_r}\sum_{\substack{\mathbf{n}=(n_1,\ldots,n_r)\in\Z_{>0}^r\\|\mathbf{n}|=|\mathbf{s}|}}C_\sigma(\mathbf{s};\mathbf{n})\polylog^\shuffle_\bullet(X_0^{n_1-1}X_{z_{\sigma(1)}}\cdots X_0^{n_r-1}X_{z_{\sigma(r)}})\notag\\
            &=\sum_{\sigma\in\mathfrak{S}_r}\sum_{\substack{\mathbf{n}\in\Z_{>0}^r\\|\mathbf{n}|=|\mathbf{s}|}}C_\sigma(\mathbf{s};\mathbf{n})\polylog^\shuffle_{\mathbf{n}}(z_{\sigma(1)},\ldots,z_{\sigma(r)})\label{eq.product(2)}.
        \end{align}
        By applying $\prod_{i=1}^r D_i^{\circ j_i}$ to both sides of the equation \eqref{eq.product(1)}=\eqref{eq.product(2)}, we have
        \begin{align}
            \prod_{i=1}^r\polylog_{s_i-j_i}(z_i)&=\left(\prod_{i=1}^rD_i^{\circ j_i}\right)\left(\prod_{i=1}^r\polylog_{s_i}(z_i)\right)\notag\\
            &=\sum_{\sigma\in\mathfrak{S}_r}\sum_{\substack{\mathbf{n}\in\Z_{>0}^r\\|\mathbf{n}|=|\mathbf{s}|}}C_\sigma(\mathbf{s};\mathbf{n})\left(\prod_{i=1}^rD_i^{\circ j_i}\right)(\polylog^\shuffle_{\mathbf{n}}(z_{\sigma(1)},\ldots,z_{\sigma(r)})).\label{eq.multi-derivation of MPL}
        \end{align}
        Here in \eqref{eq.multi-derivation of MPL}, Corollary \ref{index descending} implies
        \begin{align}
            \left(\prod_{i=1}^rD_i^{\circ j_i}\right)(\polylog^\shuffle_{\mathbf{n}}(z_{\sigma(1)},\ldots,z_{\sigma(r)}))&\overset{\text{Prop.\ref{differential operator}(i)}}{=}\left(\prod_{i=1}^rD_{\sigma(i)}^{\circ j_{\sigma(i)}}\right)(\polylog^\shuffle_{\mathbf{n}}(z_{\sigma(1)},\ldots,z_{\sigma(r)}))\notag\\
            &=\sum_{\substack{\mathbf{j}_\sigma\succeq(\mathbf{i},0),\mathbf{i}\in\Z^{r-1}_{\geq 0}}}(-1)^{\abs{\mathbf{i}}}\binom{\mathbf{j}_\sigma}{(\mathbf{i},0)}\polylog^\shuffle_{(\mathbf{i},0)-(0,\mathbf{i})-\mathbf{j}_\sigma+\mathbf{n}}(z_{\sigma(1)},\ldots,z_{\sigma(r)}).
        \end{align}
        Hence we completely prove the desired formula.
    \end{proof}
\addtocontents{toc}{\protect\setcounter{tocdepth}{0}}
    \subsection*{Relation to our previous work}
    By definition, it is clear that
    \begin{align}
        \polylog^\shuffle_{\{0\}^r}(z_1,\ldots,z_r)=\sum_{l_1,\ldots,l_r\geq 1}z_1^{l_1}\cdots z_r^{l_r}=\prod_{i=1}^r\left(\sum_{l_i\geq 1}z_i^{l_i}\right)=\prod_{i=1}^r\polylog_0(z_i).
    \end{align}
    Then the following holds:
    \begin{cor}\label{multi-variable NPMPL product}
        Let $\mathbf{j}=(j_1,\ldots,j_r)\in\Z_{\geq 0}^r$ and $z_1,\ldots,z_r\in\C$ satisfying $|z_i|<1$ for each $i$. Then we have:
        \begin{enumerate}[(1)]
            \item 
            \begin{equation}\label{eq.multi-variable NPMPL product}
                \prod_{i=1}^r\polylog_{-j_i}(z_i)=\sum_{\substack{\mathbf{j}\succeq(\mathbf{i},0),\mathbf{i}\in\Z^{r-1}_{\geq 0}}}(-1)^{\abs{\mathbf{i}}}\binom{\mathbf{j}}{(\mathbf{i},0)}\polylog^\shuffle_{(\mathbf{i},0)-(0,\mathbf{i})-\mathbf{j}}(z_1,\ldots,z_r),
            \end{equation}
            \item(\cite[Corollary 4.7]{Ki25})
            \begin{equation}\label{eq.mono-variable NPMPL product}
                \prod_{i=1}^r\polylog_{-j_i}(z)=\sum_{\substack{\mathbf{j}\succeq(\mathbf{i},0),\mathbf{i}\in\Z^{r-1}_{\geq 0}}}(-1)^{\abs{\mathbf{i}}}\binom{\mathbf{j}}{(\mathbf{i},0)}\polylog_{(\mathbf{i},0)-(0,\mathbf{i})-\mathbf{j}}(z).
            \end{equation}
        \end{enumerate}
    \end{cor}
    \begin{proof}
        (1) is clear from \eqref{eq.product of derivatives} and Corollary \ref{index descending}. (2) is obtained by setting $z_1=\cdots=z_r=:z$ in (1).
    \end{proof}
\addtocontents{toc}{\protect\setcounter{tocdepth}{1}}
    
\end{document}